\documentclass[11pt]{article}
\usepackage[utf8]{inputenc}
\usepackage{blindtext}
\usepackage{soul}
\usepackage{latexsym} % special \LaTeX symbols
\usepackage{amsmath,amssymb,amsfonts,amsthm} % AMS math
\usepackage{graphics,psfrag} % graphics support 
\usepackage{fancyhdr,lastpage} % for making nice headers
\usepackage{color,verbatim} % text markup
\usepackage[hidelinks]{hyperref} % linking and navigation
\usepackage{hyperref} % linking and navigation
\usepackage[table]{xcolor} % make colored tables
\usepackage{metalogo}
\usepackage{etoolbox}
\usepackage{cite}
\usepackage{mathrsfs}
\usepackage[affil-it]{authblk}

\usepackage{xfrac}
\usepackage{bbm}
\patchcmd{\thebibliography}{\section*{\refname}}{}{}{}

\newtheorem{theorem}{Theorem}[section]
\newtheorem{corollary}{Corollary}[theorem]
\newtheorem{lemma}[theorem]{Lemma}

\newtheorem{definition}[theorem]{Definition}
\newtheorem{example}[theorem]{Example}
\newtheorem{remark}[theorem]{Remark}

\hypersetup{
    colorlinks=false,
    pdfborder={0 0 0},
}

\usepackage{calc}

\usepackage[a4paper,top=1in,bottom=0.8in,right=0.8in,left=0.8in]{geometry}
\let\oldbibliography\thebibliography
\renewcommand{\thebibliography}[1]{%
  \oldbibliography{#1}%
  \setlength{\itemsep}{-1.5mm}%
}

\def\R{\mathbb{R}}

\def\N{\mathbb{N}}
\def\P{\mathbb{P}}

\def\Q{\mathbb{Q}}

\newcommand{\be}{\begin{equation}}
\newcommand{\ee}{\end{equation}}
\newcommand{\bea}{\begin{eqnarray}}
\newcommand{\eea}{\end{eqnarray}}
\newcommand{\beann}{\begin{eqnarray*}}
\newcommand{\eeann}{\end{eqnarray*}}
\newcommand{\benn}{\begin{equation*}}
\newcommand{\eenn}{\end{equation*}}

\newcommand{\cB}{{\mathcal B}}  % calligraphic B
\newcommand{\cC}{{\mathcal C}}  % calligraphic C
\newcommand{\cF}{{\mathcal F}}  % calligraphic F
\newcommand{\cK}{{\mathcal K}}  % calligraphic K
\newcommand{\cM}{{\mathcal M}}  % calligraphic M
\newcommand{\cP}{{\mathcal P}}  % calligraphic P
\newcommand{\cU}{{\mathcal U}}  % calligraphic U
\newcommand{\cW}{{\mathcal W}}  % calligraphic W
\newcommand{\cX}{{\mathcal X}}  % calligraphic X
\newcommand{\ie}{i.e.\ }

\newcommand{\Space}{\mathbf{Z}}
\newcommand{\Topo}{\mathcal{O}_\Space}
\newcommand{\Borel}{\mathcal{B}(\Space)}

\newcommand{\signed}{\pm}

\newcommand{\Proba}{\mathcal{M}_1(\Space)}
\newcommand{\SubProba}{\mathcal{M}_{\leq1}(\Space)}
\newcommand{\Meas}{\mathcal{M}_+(\Space)}
\newcommand{\SignedMeas}{\mathcal{M}_{\pm}(\Space)}

\newcommand{\CbFunct}{C_b(\Space)}

\newcommand{\Graphon}{\mathcal{W}_1}

\newcommand{\SubGraphon}{\mathcal{W}_{\leq 1}}
\newcommand{\Kernel}{\mathcal{W}_{\signed}}
\newcommand{\Kernelp}{\mathcal{W}_{+}}

\newcommand{\UGraphon}{\widetilde{\mathcal{W}}_1}

\newcommand{\UGraphond}{\widetilde{\mathcal{W}}_{1,d}}

\newcommand{\UKernel}{\widetilde{\mathcal{W}}_{\signed}}

\newcommand{\UKerneld}{\widetilde{\mathcal{W}}_{\signed,d}}
\newcommand{\UKernelpd}{\widetilde{\mathcal{W}}_{+,d}}
\newcommand{\UKernelp}{\widetilde{\mathcal{W}}_{+}}

\newcommand{\cKd}{\widetilde{\cK}_{d}}

\newcommand{\TotalMass}[1]{\norm{#1}_\infty}

\newcommand{\InvRelabel}{S_{[0,1]}}
\newcommand{\Relabel}{\bar{S}_{[0,1]}}

\newcommand{\NcutR}[1]{\Vert#1\Vert_{\square,\R}}

\newcommand{\NcutRSymbol}{\Vert\cdot\Vert_{\square,\R}}

\newcommand{\simd}{\sim_{d}}

\newcommand{\dd}{\delta_{\square}}

\newcommand{\NmeasF}[1]{\Vert#1\Vert_{\mathcal{F}}}

\newcommand{\NmeasFSymbol}{\Vert\cdot\Vert_{\mathcal{F}}}

\newcommand{\F}{\mathcal{F}}

\renewcommand{\P}{\mathbb{P}}

\newcommand{\un}{\mathbbm{1}}

\newcommand{\drv}{\mathrm{d}}
\newcommand{\rd}{\mathrm{d}}

\newcommand{\norm}[1]{\Vert#1\Vert}

\fancypagestyle{plain}{%
  \fancyhf{}%
  \fancyfoot[C]{\footnotesize Page \thepage\ of \pageref{LastPage}}%
}

\title{Probability graphons: the right convergence point of view}
\author{Giulio Zucal\thanks{giulio.zucal@mis.mpg.de}}
\affil[1]{Max Planck Institute for Mathematics in the Sciences, Leipzig, Germany}

\date{\today}

\begin{document}
\maketitle

\begin{abstract}
We extend the theory of probability graphons, continuum representations of edge-decorated graphs arising in graph limits theory, to the `right convergence’ point of view. First of all, we generalise the notions of overlay functionals and quotient sets to the case of probability graphons. Furthermore, we characterise the convergence of probability graphons in terms of these global quantities. In particular, we show the equivalence of these two notions of convergence with the unlabelled cut-metric convergence (and thus also with the homomorphism densities convergence and the subgraph sampling convergence). In other words, we prove the equivalence of the `left convergence’ and the `right convergence’ views on probability graphons convergence, generalising the corresponding result for (real-valued) graphons (the classical continuum representation for simple graphs).

    \vspace{0.2cm}
\noindent {\bf Keywords:}  Graph limits, large networks, probability graphons, edge-decorated graphs, dense weighted graph sequences, random matrices   

    \vspace{0.2cm}
\noindent {\bf  Mathematics Subject Classification Number:}  05C80 (Random Graphs)   60B20 (Random matrices) 60B10 (Convergence of measures)
\end{abstract}

\section{Introduction}

\subsection{Motivation and background}
Networks play a key role in the modelling, analysis, and design of complex systems in a wide range of scientific, social and technological domains, including neurobiology \cite{fornito_2016_fundamentals}, economics \cite{hausmann_2013_atlas, battiston_2012_debtrank}, urban systems \cite{barthelemy_2016_structure}, epidemiology \cite{pastor-satorras_2015_epidemic} and electrical power grids \cite{pagani_2013_power}. Graphs provide a unifying mathematical framework that allows the transfer of ideas and tools across different disciplines. However, in all these fields, the graphs under consideration are typically very large, making it impossible to use combinatorial techniques to study these networks. Moreover, for very large networks such as the Internet or the brain, information such as the exact number of nodes and other specific features of the underlying graph is not available. To deal with this complexity, graph limit theory, the study of graph sequences, their convergence and their limit objects \cite{LovaszGraphLimits}, has emerged as a powerful alternative. To name a few recent applications, graph limits have been successfully applied to nonparametric statistics \cite{wolfe2013nonparametricgraphonestimation,BorgsNonparaStat}, network dynamics \cite{Kuehn_2020GraphlimitDynam1,Kuehn_2019GraphlimitDynam2,bramburger2023pattern,MedvedevGraphLimitDynamics2,MedvedevGraphLimitDynamics1,gkogkas2022mean,bick2024dynamicalGraphLimi,ayi2023graphlimitinteractingparticle,ayi2024meanfieldlimitnonexchangeablemultiagent,kuehn2022vlasov,jabin2021meanfield} and Graph Neural Networks (GNNs) \cite{keriven2024functions,NEURIPS2023_8154c89c,maskey2022generalization,maskey2023transferability,le2023poincar,levie2024graphon}.

In graph limits theory, the limit objects are analytic structures that contain only the relevant information for a large graph. Graphons \cite{BORGS20081801,Lovsz2007SzemerdisLF, LOVASZ2006933,borgs2011convergentAnnals} are continuum representations of simple graphs used in the graph limit theory for dense simple graph sequences, i.e.\ when the simple graphs considered in the sequence have asymptotically a number of edges proportional to the square of the number of vertices. Convergence notions for graph sequences with a uniform bound on the degrees have also been studied extensively. The most famous notions of convergence for uniformly bounded degree graph sequences are Benjamini-Schramm convergence, introduced by I.\ Benjamini and O.\ Schramm \cite{BenjaminiLimit}, and its stronger version called local-global convergence \cite{local-global1,Hatami2014LimitsOL}. The reader is also referred to the monograph \cite{LovaszGraphLimits}.

In recent years, graph limits theory has developed in several directions. One important direction is the study of graph sequences of intermediate density (that are neither dense nor of uniformly bounded degree) \cite{BORGS20081801,Lovsz2007SzemerdisLF, LOVASZ2006933,borgs2011convergentAnnals,MarkovSpaces,KUNSZENTIKOVACS20191,Kunszenti_Kov_cs_2022,frenkel2018convergence,backhausz2018action}. In particular, the notion of action convergence introduced in \cite{backhausz2018action} has the advantage of bringing graphons convergence, local-global convergence and the convergence of intermediate density sequences of simple graphs into the same theoretical framework. This limit theory is based on a `right convergence'/global viewpoint and is based on functional analysis and measure theory. Other works related to this convergence notion are \cite{ArankaAction2022,MeasTheorActionZucal,zucal2023action}.

Another interesting direction is the extension of these limit notions to more general combinatorial structures such as hypergraphs and simplicial complexes \cite{hypergrELEK20121731,HypergraphsSzegedy2,HypergraphonsZhao,zucal2023action}. There has been a lot of recent interest in applications in higher order interactions (interactions beyond pairwise) and the phenomena they generate \cite{HypergraphsNetworks,HigerOrdIntBook,carletti2020dynamical,majhi2022dynamics,MHJ,HypergraphsDynamics,Bohle_2021,JOST2019870,JostMulasBook}. Hypergraphs and simplicial complexes are the natural combinatorial structures for representing higher-order interactions.

The connection between graph limit theory and the theory of random matrices is also promising \cite{RandomMatricesGraphonsZhu,RandomMatrixGraphonmale2014,backhausz2018action}.

Another direction related to random matrices, that has attracted considerable interest is the development of a limit theory for dense graphs with a decoration on the edges \cite{lovász2010limits, falgasravry2016multicolour, KUNSZENTIKOVACS2022109284,abraham2023probabilitygraphons} 
and their limit representations called probability graphons. With an edge decorated graph we mean a graph that carries some additional information on the edges, such as a weight or a colour for example. In particular, networks in applications often carry some important additional features about whether two vertices are interacting, such as how intense the interaction is. This is also the direction we focus on in this paper. In the previous works on probability graphons, the focus has been on the development of the metric, the `left convergence' and the sampling point of view, generalising the results for real-valued graphons obtained in \cite{BORGS20081801}. 
In this work, on the other hand, we develop the `right convergence' point of view extending the results in \cite{borgs2011convergentAnnals} to the case of probability graphons. In particular, we characterise the convergence of probability graphons in terms of global quantities as overlay functionals and quotient graphs.  Hence, we prove the equivalence of the `left-convergence'/local point of view established in \cite{abraham2023probabilitygraphons} and our `right convergence'/ global point of view. We explain our main result in more detail in the next section.

As already underlined, the `right convergence' point of view plays a key role in the development of sparse simple graph limits \cite{backhausz2018action, KUNSZENTIKOVACS20191}, hypergraph limits \cite{HypergraphonsZhao,zucal2023action} and recently also in the emerging matroid limits theory \cite{SubModLovaszquotientconvergence}. 

In the companion paper \cite{NewPreprintGiulio}, we use the `right convergence' viewpoint developed here to connect probability graphons convergence with action convergence \cite{backhausz2018action} and the Aldous-Hoover theorem for infinite exchangable arrays\cite{aldous1981representations,aldous2010exchangeability,hoover1979relations}.

\subsection{Main result}
Let $\Space$ be a Polish space and $[0,1]$ the unit interval endowed with the Lebesgue measure $\lambda$. A probability graphon is a measurable map from the unit square $[0,1]\times [0,1]$ (equipped with the Lebesgue measure) to the space of probability measures on $\Space$ (equipped the topology induced by the weak convergence of measures). In \cite{abraham2023probabilitygraphons}, the authors defined the convergence of probability graphons with respect to unlabelled cut-metrics $\delta_{\square}$ (with respect to different probability metrics see Definition \ref{def:ddcut} and Definition \ref{def:ddcut_cF} for example). Moreover, it was shown that this convergence is equivalent to the convergence of all homomorphism densities and sampled subgraphs convergence, see Theorem 1.5 in \cite{abraham2023probabilitygraphons}.

In this work, we generalise the notions of overlay functionals and quotient sets from the case of real-valued graphons to probability graphons. The interested reader will find the definitions, together with the combinatorial and statistical physics interpretations, of overlay functionals and quotient sets for real-valued graphons in Section 12 in \cite{LovaszGraphLimits} and in \cite{borgs2011convergentAnnals}. 
We then study the convergence of these global quantities characterizing the convergence of probability graphons. In particular, we prove in our main result that, for a sequence of probability graphons, overlay functionals convergence and quotient sets convergence are equivalent to convergence in unlabelled cut-metric (and therefore to homomorphism densities and sampled subgraphs convergence).
We briefly give more information about overlay functionals and quotient sets and our main result here. 

\textbf{Overlay functional:}

Let $\CbFunct$ be the space of real-valued continuous and bounded functions defined on $\Space.$ A $\CbFunct-$graph is a triple $G^{\beta}=(G,\beta(G))=(V(G),E(G),\beta(G))$ where $G$ is a simple graph with vertex set $V=V(G)$ and edge set $E(G)$ and $\beta(G)$ is a function
$$
\beta(G):V(G)\times V(G)\rightarrow \CbFunct $$
such that $\beta(G)_{v,w}\neq 0$ if and only if $\{v,w\}\in E(G).$

Let $[0,1]$ be the unit interval endowed with the Lebesgue measure $\lambda.$ For a probability measure $\alpha$ on $[k]$, we denote by $\Pi(\alpha)$ the set of partitions $\left\{S_1, \ldots, S_k\right\}$ of $[0,1]$ into $k$ measurable subsets with $\lambda\left(S_i\right)=\alpha_i$. For each probability graphon  $W$ and $\CbFunct-$graph $G^{\beta}$ on the vertex set $[k]$, we define the overlay functional
\begin{equation}  
\mathcal{C}(W, G^{\beta})=\sup _{\left(S_{ 1}, \ldots, S_k\right) \in \Pi(\alpha)} \sum_{i, j \in[k]}  \int_{S_i \times S_j}\int_{\Space} \beta_{i j}(G)(z)U(x, y,dz) \mathrm{d} x \mathrm{d} y,
\end{equation}
where we abbreviated with $\mathrm{d} x$ the differential $\lambda(\mathrm{d} x ).$

\textbf{Quotient graphs:}

Let $W$ be a probability graphon and let $\mathcal{P}=\left\{S_1, \ldots, S_k\right\}$ be a measurable $k$-partition of $[0,1]$ for $k\geq 1$, we define the quotient graph (or simply quotient) $W / \mathcal{P}$  as the decorated weighted graph on $[k]$, with node weights $\alpha_i(W / \mathcal{P})=\lambda\left(S_i\right)$ and as decoration of the edge $e=\{i,j\}$ the measure
$$
\beta_{i j}(W / \mathcal{P})=\frac{1}{\lambda\left(S_i\right) \lambda\left(S_j\right)} \int_{S_i \times S_j} W .
$$
We denote by $\mathcal{Q}_k(W)$ the set of quotients (or quotient set) for all $k$-partitions and we consider these as subsets of the space $$\R^k\times\cP(\Space)^{{k} \choose {2}}.$$

\textbf{Main theorem:
}
We can finally state a simplified version of the main result of this work, see Theorem \ref{ThmEquivalenceConvQuotientOverlay} and \ref{ThmEquivalenceConvQuotientOverlayVersion2} for the statement in full generality.

\begin{theorem}
For any sequence of probability graphons $\left(W_n\right)$ and a probability graphon $W$, the following are equivalent:
\begin{enumerate}

\item  The sequence $\left(W_n\right)$ is convergent to $W$ in the cut distance $\delta_{\square};$
\item The overlay functional values $\mathcal{C}\left(W_n, G^{\beta}\right)$ converge to $\mathcal{C}\left(W, G\right)$ for every $\CbFunct-$graph $G^{\beta}$;
\item The quotient sets $\mathcal{Q}_k\left(W_n\right)$ converge to $\mathcal{Q}_k\left(W\right)$ in the $d_{\square}^{\text {Haus }}$ Hausdorff metric (see  equation \eqref{Eqdsquare} and Definition \ref{DefHausdorffDist}) for every $k \geq 1$.
\end{enumerate}
\end{theorem}

Therefore, in this work, we show the equivalence of convergence in unlabelled cut-metric, convergence of overlay functionals and convergence of quotient sets. Moreover, our result together with Theorem 1.5 in \cite{abraham2023probabilitygraphons} shows that the `right convergence'/global point of view developed in this work is equivalent to the `left convergence'/local point of view developed in \cite{abraham2023probabilitygraphons} and \cite{lovász2010limits}.

We note that in the statements in full generality, Theorem \ref{ThmEquivalenceConvQuotientOverlay} and Theorem \ref{ThmEquivalenceConvQuotientOverlayVersion2}, we also characterize the convergence of overlay functionals in a less wasteful way, using only $\CbFunct-$graphs $G^{\beta}$    

Our proofs build on results from \cite{abraham2023probabilitygraphons} and follow the non-effective proof scheme developed for real-valued graphons from \cite{borgs2011convergentAnnals} which is also nicely explained in Chapter 12.3 in \cite{LovaszGraphLimits}. To apply this proof scheme we overcome many difficulties caused by the infinite-dimensional and metric nature of the probability graphons setting, as opposed to the one-dimensional space $\R$ equipped with the Euclidean norm considered in the real-valued graphons setting.
\newline

\subsection{Organization}
In Section \ref{Sec2Notat}, the notation is fixed and basic measure-theoretic notions needed in this work are introduced. In Section \ref{Sec3ProbGraphon}, the theory of probability graphons is summarised. Finally, in Section \ref{Sec4Results}, we introduce new invariants for probability graphons, overlay functionals and quotient sets, study their properties and prove the main result of this work, a characterisation of the convergence of probability graphons from the global/ `right convergence' point of view.  
\section{Notation and basic notions}\label{Sec2Notat}
This section introduces the basic mathematical concepts used in this work and defines the notation.
\subsection{Measure theory background}

In this section, we introduce basic measure-theoretic concepts and fix the notation following \cite{abraham2023probabilitygraphons}.

A function $\varphi : \Omega_1 \to \Omega_2$ between two probability spaces $(\Omega_i, \mathcal{A}_i, \P_i)$, $i = 1, 2$, is called measure-preserving if it is measurable and if for every set $A \in \mathcal{F}_2$ the condition $\P_2(A) = \P_1(\varphi^{-1}(A))$ holds. This is equivalent to saying that for any measurable non-negative function $f : \Omega_2 \to \mathbb{R}$, the following equation is satisfied:
\[
\int_{\Omega_1} f(\varphi(x))\, d\P_1(x) = \int_{\Omega_2} f(x)\, d\P_2(x).
\]
In this work, we will always consider the unit interval $[0,1]$ and the unit square $[0,1]^2=[0,1]\times [0,1]$ equipped with the Euclidean norm, the Borel $\sigma-$algebra and the Lebesgue measure. These are probability spaces.
We denote by $\InvRelabel$ the set of all bijective measure-preserving maps from $[0,1]$ to itself, while $\Relabel$ denotes the set of measure-preserving maps from $[0,1]$ to itself.
\medskip

Let $(\Space,d_{\Space})$ be a non-empty Polish metric space. Let $\Topo$ be the topology generated by the metric $d_{\Space}$ and $\Borel$ be the Borel $\sigma$-algebra on $\Space$ generated by the topology $\Topo$. Furthermore, we let $\CbFunct$ be the  space of continuous bounded real-valued functions on $(\Space,\Topo)$ equipped with the supremum norm $\|\cdot\|_{\infty}.$ We denote by $\SignedMeas$ the space of finite signed measures on $(\Space,\Borel)$, $\Meas$ as the subspace of measures, $\SubProba$ as the subspace of measures with total mass at most $1$, and $\Proba$ as the subspace of probability measures. The following chain of inclusions holds:
\[
  \Proba \subset \SubProba \subset \Meas \subset \SignedMeas .
\]
For a signed measure $\mu\in\SignedMeas$, we recall the unique Hahn-Jordan decomposition $\mu  = \mu^+ -  \mu^-$ where $\mu^+, \mu^-  \in\Meas$ are mutually singular measures (that  is $\mu^+(A)=0$ and $\mu^-(A^c)=0$ for some measurable set $A$). The total variation measure $\vert\mu\vert$ of $\mu\in\SignedMeas$ is defined as $\vert\mu\vert  =   \mu^+  +  \mu^-  \in  \Meas$. In particular, for a measure $\mu\in\Meas$, we have $\vert  \mu \vert  = \mu$. For a signed measure $\mu\in \SignedMeas$ and a real-valued measurable function $f$ defined on $\Space$, we denote by $\mu[f]=\mu(f)=\langle \mu, f \rangle=
\int f \, \rd \mu=\int_{\Space} f(x)\, \rd\mu$ the integral of $f$ with respect to $\mu$ when well defined. For a signed measure  $\mu\in\SignedMeas$, the total mass of $\mu$ is $\|{\mu}\|_{TV}=\mu^+(\Space) + \mu^-(\Space),$ which also equals the supremum of $\mu(f)$ taken over all measurable functions $f$ with values in $[-1,1]$.

We consider the space $\SignedMeas$ equipped with the weak convergence topology, i.e. the smallest topology such that for each $f\in\CbFunct$ the map $\mu \mapsto \mu(f)$ is continuous. In particular, a sequence  of signed measures $(\mu_n)_{n\in\N}$  weakly converges to some $\mu\in\SignedMeas$ if and only if, for each function $f\in   \CbFunct$, we have $\lim_{n\to +\infty} \mu_n(f) = \mu(f)$. Recall that $\Meas$ and $\Proba$ equipped with the weak convergence topology are Polish spaces. 

The maps $\mu   \mapsto  \mu^+$ and $\mu   \mapsto  \mu^-$ (and thus also $\mu \mapsto \vert\mu\vert$) are measurable (see \cite[Theorem~2.8]{MeasurableSetsMeasures} and Remark 2.4 in \cite{abraham2023probabilitygraphons}). As a consequence, the map $\mu \mapsto \|\mu\|_{TV}$ is also measurable.
Observe that $\Proba$ and $\Meas$ are closed, and therefore measurable subsets of $\SignedMeas.$

\begin{definition}
    A sequence of $[0,1]$-valued functions $\F = (f_k )_{k\in\N}$ in $\CbFunct$, with $f_0=\un_{\mathbf{Z}}$ the constant function equal to one, is    
\begin{enumerate}
\item  \textbf{Separating}   if  for all measures $\mu,\nu$ from $\SignedMeas$ (or equivalently just from  $\Meas$) such that for every $k\in\N$, $\mu(f_k) = \nu(f_k)$, then $\mu = \nu$.
 \item  \textbf{Convergence determining}  if   for  every $(\mu_n)_{n\in\N}$  and  $\mu$   measures  in  $\Meas$  such  that we have $\lim_{n\to +\infty} \mu_n(f_k)  = \mu(f_k)$ for all $k\in\N$, then $(\mu_n )_{n\in\N}$ weakly converges to $\mu$.
\end{enumerate}
\end{definition}
On Polish spaces, there always exists a convergence determining sequence, see \cite[Corollary 2.2.6]{Bogachev} or the proof of Proposition 3.4.4 in \cite{Ethier} and  Remark 2.3 in \cite{abraham2023probabilitygraphons}. Moreover, observe that convergence determining sequences are separating. 

The interested reader can find more details about the weak convergence of signed measures in Section 2 in \cite{abraham2023probabilitygraphons} and in the standard references \cite{Bogachev,BogachevMT1,BogachevMT2,Ethier,MeasurableSetsMeasures,Varadarajan}.

\subsection{A metric metrizing weak convergence of measures:}

 In this subsection, we recall the notion of Lévy-Prokhorov Metric, a metric on the space of measures with the important property of metrizing weak convergence of measures. 
 
\begin{definition}[Lévy-Prokhorov metric]\label{LevyProk}
 The \emph{Lévy-Prokhorov Metric} $d_{\mathrm{LP}}$ on the space of measures $\Meas$ is for $\eta_1,\eta_2\in \Meas$
$$\begin{aligned}
d_{\mathrm{LP}}\left(\eta_{1}, \eta_{2}\right)=&\inf \left\{\varepsilon>0: \eta_{1}(U) \leq \eta_{2}\left(U^{\varepsilon}\right)+\varepsilon \text{ and } \right.\\
&\left.\eta_{2}(U) \leq \eta_{1}\left(U^{\varepsilon}\right)+\varepsilon  \text{ for all } U \in \Topo\right\},
\end{aligned}$$

where we recall that $\Topo$ is the Borel $\sigma$-algebra on $\Space$ generated by $d_{\Space}$, $U^{\varepsilon}$ is the set of points that have a distance $d_{\Space}$ smaller than $\varepsilon$ from $U$.
\end{definition}

In particular, a sequence $(\mu_n)_{n\in \N}$ in $\Meas$ is weak convergent if and only if it is convergent in Levy-Prokhorov metric, see Theorem 6.8 in \cite{billingsley1968convergence}. Moreover, the Levy-Prokhorov metric is upper-bounded by the total variation.
\newline
\begin{lemma}[Page 12 in \cite{ComparingProbMeas}] \label{LemmBoundProhTotVar}
For any two measures $\mu,\nu \in \mathcal{M}_+$, we  have  that Levy-Prokhorov metric is bounded by the total variation distance in the following way $$d_{\mathcal{LP}}(\mu,\nu)\leq \|\mu-\nu\|_{TV}.$$
\end{lemma}

Therefore, if we compare only probability measures, the Lévy-Prokhorov metric is upper-bounded by $1.$

We now prove a simple lemma which will be useful in the next sections.

\begin{lemma}\label{LemmaIneqScalingProkhorov}
  Let $(X,d)$ be a metric space. Let $\mu,\nu$ two measures on $X$ and $\alpha >1$. Then 
    $$
    d_{\mathcal{LP}}(\mu,\nu)\leq d_{\mathcal{LP}}(\alpha\mu,\alpha\nu)\leq \alpha d_{\mathcal{LP}}(\mu,\nu).$$ Moreover, the previous inequalities are sharp.
\end{lemma}
\proof
We first prove $d_{\mathcal{LP}}(\alpha\mu,\alpha\nu)\leq \alpha d_{\mathcal{LP}}(\mu,\nu)$. Let $\varepsilon >d_{\mathcal{LP}}(\mu,\nu) ,$ then $\mu(B)\leq \nu (B^{\varepsilon})+\varepsilon$ by the definition of the Levy-Prokhorov distance \ref{LevyProk}. Therefore, $$
\alpha \mu(B)\leq \alpha (\nu(B^{\varepsilon})+\varepsilon)=
\alpha\nu(B^{\varepsilon})+ \alpha\varepsilon \leq \alpha\nu(B^{\alpha\varepsilon})+ \alpha\varepsilon. 
$$ 
Sending $\varepsilon\rightarrow d_{\mathcal{LP}}(\mu,\nu),$ as we can also exchange the role of $\mu $ and $\nu$, we obtain the inequality. Consider now $X=\R$ with the Euclidean distance and $\mu=\delta_0$ and $\mu=\delta_{10\alpha}$. We obtain $d_{\mathcal{LP}}(\mu,\nu)=1$ and $d_{\mathcal{LP}}(\alpha\mu,\alpha\nu)=\alpha=\alpha d_{\mathcal{LP}}(\mu,\nu)$ for any $\alpha>1.$ This proves that the inequality is sharp.

We now show that $d_{\mathcal{LP}}(\mu,\nu)\leq d_{\mathcal{LP}}(\alpha\mu,\alpha\nu)$.  Let $\varepsilon >d_{\mathcal{LP}}(\alpha\mu,\alpha\nu) ,$ then $\alpha\mu(B)\leq \alpha\nu (B^{\varepsilon})+\varepsilon$ by definition of Levy-Prokhorov distance \ref{LevyProk}. Therefore, $$
\mu(B)\leq  \nu(B^{\varepsilon})+\frac{\varepsilon}{\alpha} \leq \nu(B^{\varepsilon})+ \varepsilon. 
$$ 
Sending $\varepsilon\rightarrow d_{\mathcal{LP}}(\alpha\mu,\alpha\nu),$ as we can also exchange the role of $\mu $ and $\nu$, we obtain the inequality. Consider now $X=\R$ and $\mu=\delta_0$ and $\nu=\delta_{1}.$ We find that $d_{\mathcal{LP}}(\mu,\nu)=1$ and $d_{\mathcal{LP}}(\alpha\mu,\alpha\nu)=1= d_{\mathcal{LP}}(\mu,\nu)$ for any $\alpha>1.$ Therefore, this inequality is also sharp.
\endproof

Another useful property of the Levy-Prokhorov is given in the following lemma.

%{\color{blue} Lemma 3.21 in \cite{abraham2023probabilitygraphons}}

%{\color{blue} \cite[Lemma~3.21]{abraham2023probabilitygraphons}}

\begin{lemma}[Lemma 3.21 in \cite{abraham2023probabilitygraphons}]\label{LemmaQuasi-convProhorov}
    The Levy-Prokhorov distance $d_{\mathcal{LP}}$ is quasi-convex on $\mathcal{M}_{+}$, i.e.\ for any measures $\mu_1,\mu_2,\nu_1,\nu_2 \in \mathcal{M}_{+} $ and any $\alpha \in [0,1]$ we have $$
    d_{\mathcal{LP}}(\alpha \mu_1+(1-\alpha)\mu_2,\alpha \nu_1+(1-\alpha)\nu_2)\leq \max(d_{\mathcal{LP}}(\mu_1, \nu_1),d_{\mathcal{LP}}( \mu_2,\nu_2)).    
    $$
\end{lemma}

\subsection{A norm  based on a convergence determining sequence}
	\label{subsection_norm_F}

In this section, we present another important metrization of weak convergence on the space of measures.

Let $\F=(f_k  )_{k\in\N}$ in $\Meas$ be a convergence determining sequence. We recall that we required in the definition of convergence determinig sequence that $f_0=\mathbbm{1}_{[0,1]}$ and $f_k \in\CbFunct$ takes values in $[0, 1]$ for each $k\in\N.$ For $\mu\in \SignedMeas$, the norm $\|\cdot\|_{\cF}$ generated by the convergence determining sequence $\cF$ is defined as follows: 
\begin{equation}\label{eq_def_NmeasF}
\NmeasF{\mu}\ = \sum_{k\in \N} 2^{-k} |\mu(f_k)|.
\end{equation}
Observe that from \eqref{eq_def_NmeasF}, it follows $\NmeasF{\mu} \leq 2\|\mu\|_{TV}.$ The metric $d_{\cF},$ generated by the norm $\|\cdot\|_{\cF}$, metrizes the weak topology on $\Meas$ (not on all of $\SignedMeas$).

The norm $\NmeasFSymbol$  is not complete when  $\Space$ is
not compact (see Lemma  3.23 in \cite{abraham2023probabilitygraphons}). However, the norm $\NmeasFSymbol$ and the associated distance $d_{\cF}$ play a key role in the theory of probability graphons, see \cite{abraham2023probabilitygraphons} and the rest of this work.

\section{Probability graphons/decorated graph limits}\label{Sec3ProbGraphon}

We summarize in this section the theory of probability graphons and decorated graph limits following \cite{abraham2023probabilitygraphons}. 
\subsection{Definition of probability graphons}\label{section_def_graphons}

We briefly recall here the notion of (real-valued) graphon from dense graph limit theory \cite{LovaszGraphLimits}.
\begin{definition}[Real-valued graphon]
    A \emph{real-valued kernel} is a measurable function
    $$w:[0,1]\times [0,1]\rightarrow \R.$$
    A \emph{real-valued graphon} is a measurable function
    $$w:[0,1]\times [0,1]\rightarrow \R,$$
    i.e.\ a real-valued kernel taking values in $[0,1].$ 
\end{definition}
\begin{remark}
Sometimes, one requires real-valued graphons to be symmetric in the two variables. We do not impose this condition here. 
\end{remark}

In this work, we are interested in generalizations of real-valued kernels and graphons: measure-valued kernels and probability graphons.
Recall that $\Space$ is a Polish metric space and $\SignedMeas$ is the space of finite signed measures on $\Space$. In this work, we will focus our attention on the following mathematical objects.

\begin{definition}[Signed measure-valued kernels]\label{Def:MeasValuedKernel}
A \emph{signed measure-valued kernel} or \emph{$\SignedMeas$-valued kernel}
is a map $W$ from $[0,1]^2$  to $\SignedMeas$, 
such that:
\begin{enumerate}
\item   $W$  is   a  \emph{signed-measure}   in  $\rd   z$, i.e.\  for   every
  $(x,y) \in [0,1]^2$, $W(x,y;\cdot)$ belongs to $\SignedMeas$.
  
\item $W$ is \emph{measurable} in $(x,y)$, i.e.\ for every measurable set $A\subset \Space$, the function $(x,y)\mapsto W(x,y;A)$ defined on $[0,1]^2$ is measurable.
\item $W$ is \emph{bounded}:
\begin{equation}
   \label{eq:def:TM}
\sup\|W\|_{TV}:=    \sup_{x,y\in [0, 1] }\, \|W(x, y; \cdot)\|_{TV} <+\infty,
  \end{equation}  
  where with  $\sup_{x,y\in [0, 1] }$ we mean the essential supremum on $[0,1]^2$ with respect to the Lebesgue measure.
\end{enumerate}
\end{definition}
We will be particularly interested in the case of probability measure-valued  kernels, i.e. in signed measure-valued kernels taking values in the probability space $\Proba.$ We will call these objects \emph{probabilty graphons}. We denote by  $\Graphon$ (respectively $\SubGraphon$, respectively $\Kernelp$, respectively   $\Kernel$) the space of probability graphons (respectively sub-probability measure-valued kernels, respectively  measure-valued  kernels, respectively signed measure-valued kernels), where we identify kernels that are equal almost everywhere on $[0,1]^2$ (with respect to the Lebesgue measure). 

For  $\cM\subset\SignedMeas$, we denote by $\cW  _\cM$ the subset of signed measure-valued kernel $W\in  \Kernel$ which are $\cM$-valued: $W(x,y; \cdot)\in \cM$ for every $(x,y)\in [0, 1]^2.$

\begin{example}[On real-valued kernels]\label{rem:real-valued-kernels}
Every real-valued  graphon  $w$ can be represented as a probability graphon $W$ in the following way. Let's consider $\Space =\{ 0, 1\}$ with the discrete metric and the probability-graphon $W$ defined as $W(x,y;\drv z)  = w(x,y)  \delta_1(\drv z) +  (1-w(x,y)) \delta_0(\drv z)$ for every $x,y\in[0,1]$, where $\delta_z$ is the Dirac mass located at $z$. In particular, we have $$w(x,y)=W(x,y; \{1\})=\int_{\Space}W(x,y,\drv z)$$ for $x,y \in [0, 1]$. 
\end{example}

For a signed measure-valued kernel $W \in \Kernel,$ we denote by $W^+$ the map $W^+ : [0,1]^2 \to \Meas$ to be the positive part of $W,$ i.e.\ for every $(x,y)\in [0,1]^2$, $W^+(x,y;\cdot)$ is the positive part of the measure $W(x,y;\cdot)$. Similarly we denote $W^- : [0,1]^2 \to \Meas$ the negative part of $W.$
Moreover, we define $\vert W \vert = W^+ + W^-$ the variation of $W$ and $\| W \|_{TV} = \vert W \vert (\Space)$ the total mass/variation of $W.$

%{\color{blue} Lemma 3.3 in \cite{abraham2023probabilitygraphons}.}

%{\color{blue} \cite[Lemma~3.3]{abraham2023probabilitygraphons}}
\begin{lemma}[Lemma 3.3 in \cite{abraham2023probabilitygraphons}]\label{lem:mesurability_W_+}
The maps $W^+$, $W^-$ and $\vert W \vert$ are measure-valued kernels, and the map $\Vert W \Vert_{TV} : (x,y) \mapsto\|W(x,y;\cdot)\|_{TV}$ is measurable.
\end{lemma}

We say that a subset  $\cK\subset \Kernel$ is \emph{uniformly bounded} if:
\begin{equation}\label{eq_def_TM_graphon}
\sup_{W\in \cK}\,\, \sup \|W\|_{TV}  < +\infty .
\end{equation}

As an example, $\Graphon$ is uniformly bounded with $\sup_{W\in \Graphon}\,\, \sup\|W\|_{TV}  =1.$

\begin{remark}[Probability-graphons $W : \Omega\times \Omega \to \Proba$]
	\label{rem_vertex_type_omega}
Similarly to the case of real-valued graphons, one could replace the space $[0,1]$ by any standard probability space $(\Omega,\mathcal{A},\P),$ i.e.\\ a probability space where there exists a measure-preserving map $\varphi : [0,1] \to \Omega$, where $[0,1]$ is equipped with the Borel $\sigma$-algebra and the Lebesgue measure. In fact, sometimes it is convinient to consider probability-graphons of the form $W : \Omega\times \Omega \to \Proba$.

The measure-preserving map $\varphi$ can be used to construct a version $W^\varphi$ of $W$ on $\Omega' = [0,1]$. Pushing this approach forward, one can modify the definition of the unlabelled cut distance $\delta_{\square,\mathcal{LP}}$ (see Definition \ref{def:ddcut}) similarly to \cite[Theorem 6.9]{jansonGraphonsCutNorm2013} allowing probability-graphons to be constructed on different standard probability spaces. This was also observed in Remark 3.4 in \cite{abraham2023probabilitygraphons}. However, we focus for simplicity on the equivalent case where each probability-graphon is defined on $\Omega = [0,1]$ in this work.
\end{remark}

We say that a finite signed measure-valued kernel (or a probability graphon) $W$ is symmetric if for almost every $x,y \in [0,1]$, $$W(x,y;\cdot)=W(y,x;\cdot).$$

\begin{remark}[Symmetric kernels]
In our work, we consider general (not necessarily symmetric) measure-valued kernels and probability-graphons as in \cite{abraham2023probabilitygraphons}. However, one could consider the more restrictive case of symmetric measure-valued kernels.  
\end{remark}

We define finite signed measure-valued kernels which are stepfunctions often used for approximation.
 
\begin{definition}[Signed measure-valued stepfunctions]
  \label{def:stepfunction}
  A signed measure-valued kernel $W\in \Kernel$ is a \emph{step-function} if there exists a
  finite partition of $[0,1]$ into measurable (possibly empty) sets, say
  $\mathcal{P}=\{S_1,\cdots,S_k\}$, such that $W$ is constant on the sets
  $S_i \times  S_j$, for $1\leq i,j\leq  k$.  \end{definition}
We want to emphasise that step functions are a strict subset of simple functions. However, we will see that, for our purposes, we can use step functions interchangeably with simple functions.

\subsection{Cut distances}\label{section_def_dcut}

Following \cite{abraham2023probabilitygraphons}, we introduce a pseudo-distance and a norm for probability graphons (or more general finite signed measure-valued kernels): the \emph{cut distance} and the \emph{cut norm} respectively. These are generalizations of the cut norm for real-valued graphons and kernels, see \cite[Chapter 8]{LovaszGraphLimits}. For a finite signed measure-valued kernel $W\in\Kernel$ and a measurable subsets of the unit interval $A\subset [0,1]^2,$ we define $W(A;\cdot),$ the finite signed measure on $\Space$ given by:
$$
  W(A;\cdot) = \int_{A} W(x,y;\cdot)\ \mathrm{d} x \mathrm{d} y.
$$
 The first pseudo-distance is based on the Levy-Prokhorov metric $d_{\mathcal{LP}}.$

\begin{definition}[The cut semi-distance $d_{\square, \mathcal{LP}}$]
Let $d_{\mathcal{LP}}$ be the Levy-Pokhorov metric on $\Meas.$  %containing the zero measure. 
The associated cut semi-distance $d_{\square, \mathcal{LP}}$ is the function from $\Kernelp^2$ to $\R_{+}$ defined by:
\begin{equation}\label{cutsemi-disteq}
    d_{\square, \mathcal{LP}}(U, W)=\sup _{S, T \subset[0,1]} d_{\mathcal{LP}}(U(S \times T ; \cdot), W(S \times T ; \cdot)),
\end{equation}
where the supremum is taken over all measurable subsets $S$ and $T$ of $[0,1]$.
\end{definition}

The second pseudo-distance $d_{\square,\cF}$ (or related pseudo norm $\|\cdot\|_{\square,\cF}$) is based on the norm $\|\cdot\|_{\mathcal{F}}$ constructed using a convergence determining sequence $\cF=(f_k)_{k\in\N}$.

\begin{definition}[The cut semi-distance $d_{\square, \mathcal{F}}$]
     Let $\|\cdot\|_{\mathcal{F}}$ be the norm on $\mathcal{M}_{+/-}(\mathbf{Z})$ given by a convergence determining sequence $\cF=(f_s)_{s\in\N}$.  %containing the zero measure. 
     The associated cut semi-distance $d_{\square, \cF}$ and cut semi-norm $\|\cdot\|_{\square,\cF}$ is the function from $\Kernel$ to $\R_{+}$ defined by:
\begin{equation}\label{cutsemi-disteq_cF}
    d_{\square,\cF}(U, W)=\|U-V\|_{\square,\cF}=\sup _{S, T \subset[0,1]} \|U(S \times T ; \cdot)- W(S \times T ; \cdot))\|_{\cF},
\end{equation}
where the supremum is taken over all measurable subsets $S$ and $T$ of $[0,1]$.
\end{definition}
\medskip
For  $W\in\Kernel$ and  $f\in\CbFunct$, we denote by $W[f]$ the real-valued kernel defined by:
\begin{equation}
   \label{eq:notation_graphon_function}
  W[f](x,y) = W(x,y;f) = \int_{\Space} f(z)\ W(x,y;\drv z).
\end{equation}
We denote by $\NcutRSymbol$ the real-valued cut norm defined as:
\begin{equation}
	\label{eq_def_NcutR}
\NcutR{w} = \sup_{S,T \subset [0,1]} \left\vert \int_{S\times T} w(x,y)\ \drv x\drv y \right\vert
\end{equation}
where $w$ is a real-valued kernel. See \cite{LovaszGraphLimits} for more details about the real-valued cut-norm.

We prove the following simple Lemma relating $\|\cdot\|_{\square,\cF}$ and $\NcutRSymbol.$

\begin{lemma}\label{lemmaBoundConvDetSequenceCutNormF}

Let $(f_n)_{n\in \N}=\cF$ be a convergence determining sequence. The following inequality holds:
$$    \NcutR{W[f_n]}\leq 2^n\|W\|_{\square,\cF}.$$
\end{lemma}
\proof

From the inequality 

$$
\sum^\infty_{n=0}2^{-n}\left|\int_{S\times T}W(x,y;f_n)\ \drv x \drv y\right|\geq2^{-n}\left|\int_{S\times T}W(x,y;f_n)\ \drv x \drv y\right|$$
we get
$$\begin{aligned}
&\|W\|_{\square,\cF}=\sup_{S,T\subset[0,1]}\sum^\infty_{n=0}2^{-n}\left|\int_{S\times T}W(x,y;f_n)\ \drv x \drv y\right|\geq \\
& 2^{-n}\sup_{S,T\subset[0,1]}\left|\int_{S\times T}W(x,y;f_n)\ \drv x \drv y\right|=2^{-n}\NcutR{W[f_n]}.\end{aligned}$$

\endproof

Recall $\Relabel$ denotes the set of measure-preserving maps from $[0, 1]$ to $[0, 1]$, and $\InvRelabel$ denotes the set of bijective measure-preserving maps.

The relabeling of a signed measure-valued kernel $W$ by a measure-preserving map $\varphi\in \Relabel$, is the signed measure-valued kernel $W^\varphi$ defined for every $x,y\in [0,1]$ and every measurable set $A\subset \Space$ by:
\[
  W^\varphi(x,y;A) = W(\varphi(x),\varphi(y);A)
  \quad \text{for $x,y\in [0,1]$ and $ A\subset \Space$ measurable}.
\]
(Similarly, for a real-valued kernel $w$ the real-valued kernel $w^{\varphi}$ is defined as $w^{\varphi}(x,y)=w(\varphi(x),\varphi(y))$ for every $x,y\in [0,1].)$

We observe that $$
\sup \|W\|_{\infty}=\sup \|W^{\varphi}\|_{\infty}.$$

\subsection{Unlabeled  cut distances}
	\label{subsection_unlabeled_cut_distance}
We can now introduce the unlabelled cut distance, a premetric on the space of probability graphons (or more generally finite signed measure valued kernels).

\begin{definition}[The unlabeled cut distance $\dd$]\label{def:ddcut}
  Set  $\cK\in\{\Graphon,  \Kernelp\}$. The \emph{unlabelled cut distance} on $\cK$ is the premetric $\delta_{\square,\mathcal{LP}}$ on $\cK$ such that for $U,W \in\cK$
\begin{equation}\label{def_ddcut}
\delta_{\square,\mathcal{LP}}(U,W) = \inf_{\varphi\in\InvRelabel} d_{\square,\mathcal{LP}}(U,W^\varphi)	
= \inf_{\varphi\in\InvRelabel} d_{\square,\mathcal{LP}}\left(U^{\varphi},W\right)		 .
\end{equation}
\end{definition}

\begin{definition}[The unlabeled cut distance $\delta_{\square,\cF}$]\label{def:ddcut_cF}
Let $\cF$ be a convergence determining sequence and set $\cK\in\{\Graphon,  \Kernelp,  \Kernel\}$. The \emph{unlabelled cut distance} with respect to $\cF$ on $\cK$ is the premetric $\delta_{\square,\cF}$ on $\cK$ such that for $U,W \in\cK$ \begin{equation}\label{def_ddcut_cF}
\delta_{\square,\cF}(U,W) = \inf_{\varphi\in\InvRelabel} d_{\square,\cF}(U,W^\varphi)	
= \inf_{\varphi\in\InvRelabel} d_{\square,\cF}\left(U^{\varphi},W\right)		 .
\end{equation}
\end{definition}

Let $u,w$ two real-valued kernels one can similarly define the real-valued cut distance (see \cite{LovaszGraphLimits} for more details about the real-valued cut-distance) 
$$\delta_{\square,\R}(u,w):=\inf_{\varphi\in \InvRelabel}\NcutR{u-w^{\varphi}}.$$
Similarly to the cut distance $\|\cdot\|_{\square,\cF}$ and the real-valued cut-norm $\|\cdot\|_{\square,\R}$, we can link the metrics $\delta_{\square,\cF}$ and $\delta_{\square,\R}$ as follows. 

Let $(f_n)_{n\in \N}=\cF$ be a convergence determining sequence and $U,W\in \Kernel.$ From Lemma \ref{lemmaBoundConvDetSequenceCutNormF} it directly follows that 
$$  \delta_{\square,\R}  (W[f_n],U[f_n])\leq 2^n \delta_{\square,\cF}(U,W).$$

Observe that  $\delta_{\square,\mathcal{LP}}$ (respectively $\delta_{\square,\cF}$) is symmetric (see also Section 3.4 in \cite{abraham2023probabilitygraphons}) and satisfies the triangular inequality. Therefore, the premetric $\delta_{\square,\mathcal{LP}}$ (respectively $\delta_{\square,\cF}$) induces a distance, that we  still denote  by $\delta_{\square,\mathcal{LP}}$ (respectively $\delta_{\square,\cF}$), on  the quotient space $\cKd = \cK / \simd$ of  the space of kernels $\cK\subset \Kernel$  with respect to the equivalence relation $\simd$ given by $U\simd W$  if and only if  $\delta_{\square,\mathcal{LP}}(U,W)=\delta_{\square,\cF}(U,W)=0$.

\subsection{Weak isomorphism}
\label{section_weak_isomorphism}
We introduce the notion of weak isomorphism, an equivalence relation between probability graphons. 

\begin{definition}[Weak isomorphism]\label{def_weak_isomorphism}
We say that two signed measure-valued kernels $U$ and $W$ are \emph{weakly isomorphic}
(denoted by $U\sim W$) if there exist two measure-preserving maps $\varphi, \psi\in\Relabel$
such that $U^\varphi(x,y; \cdot) = W^\psi(x,y;\cdot)$ for almost every $x,y\in [0,1]$.

We denote by $\UKernel=\Kernel / \sim$ (resp. $\UGraphon=\Graphon / \sim$) the quotient space of signed measure-valued kernels (respectively probability-graphons) with respect to the equivalence relation $\sim$, i.e.\ the space of finite signed measure-valued kernels (respectively probability-graphons) where weakly isomorphic finite signed measure-valued kernels (respectively weakly isomorphic probability-graphons) are identified. 
\end{definition}

Observe that $U\sim W$ implies that $\sup\|U\|_{TV}=\sup\|W\|_{TV}.$
Therefore, the notion of a uniformly bounded subset from \eqref{eq_def_TM_graphon} applies naturally also to $\UKernel$.

We summarize here some additional results about weak isomorphism of probability graphons.

%{\color{blue} Theorem 3.17 in \cite{abraham2023probabilitygraphons}}

%{\color{blue} \cite[Theorem~3.17]{abraham2023probabilitygraphons}}

\begin{theorem}[Lemma 3.21 in \cite{abraham2023probabilitygraphons}]
   \label{theo:Wm=W}
   Two finite signed measures-valued kernels are weakly isomorphic, \ie $U \sim W,$ if and only if $U \simd W,$ i.e.\ when $\delta_{\square,\mathcal{LP}}(U,W)=\delta_{\square,\cF}(U,W)= 0$.

   Furthermore, $\delta_{\square,\mathcal{LP}}$ and $\delta_{\square,\cF}$ are distances on $\UGraphon=\UGraphond$ (respectively $\UKernelp=\UKernelpd$ or $\UKernel=\UKerneld$).
    %Furthermore, $\delta_{\square,\mathcal{LP}}$ and $\delta_{\square,\cF}$ are distances on $\UGraphon=\UGraphond.$ More generally, $\delta_{\square,\mathcal{LP}}$ is a distance on $\UKernelp=\UKernelpd$ and $\delta_{\square,\cF}$ is a distance on $\UKernel=\UKerneld$
\end{theorem}

The following lemma is a special case of Proposition 3.18 in \cite{abraham2023probabilitygraphons}.

%{\color{blue} Proposition 3.18 in \cite{abraham2023probabilitygraphons}}
%{\color{blue} \cite[Proposition~3.18]{abraham2023probabilitygraphons}}

\begin{lemma}[Proposition 3.18 in \cite{abraham2023probabilitygraphons}]\label{thm_min_dist}
  For the unlabelled cut distances $\delta_{\square,\mathcal{LP}}$ we have the following alternative formulations on $\Graphon$ (resp. $\Kernelp$ or $\Kernel$):
\begin{equation}
  \label{eq_premetric}
  \begin{aligned}
\delta_{\square,\mathcal{LP}}(U,W) 
& = \underset{\varphi\in \InvRelabel}{\inf} d_{\square,\mathcal{LP}}(U,W^\varphi) 
  = \underset{\varphi\in \Relabel}{\inf} d_{\square,\mathcal{LP}}(U,W^\varphi)\\
& = \underset{\psi\in \InvRelabel}{\inf} d_{\square,\mathcal{LP}}(U^\psi,W) 
= \underset{\psi\in \Relabel}{\inf} d_{\square,\mathcal{LP}}(U^\psi,W) 	\\
& = \underset{\varphi, \psi\in \InvRelabel}{\inf} d_{\square,\mathcal{LP}}(U^\psi,W^\varphi) 
= \underset{\varphi,\psi\in \Relabel}{\min} d_{\square,\mathcal{LP}}(U^\psi,W^\varphi).
\end{aligned}
\end{equation}

Similarly, for the unlabelled cut distances $\delta_{\square,\cF}$ we have the following alternative formulations on $\Graphon$ (resp. $\Kernelp$ or $\Kernel$): 
\begin{equation}
  \label{eq_premetric2}
  \begin{aligned}
\delta_{\square,\cF}(U,W) 
& = \underset{\varphi\in \InvRelabel}{\inf} d_{\square,\cF}(U,W^\varphi) 
  = \underset{\varphi\in \Relabel}{\inf} d_{\square,\cF}(U,W^\varphi)\\
& = \underset{\psi\in \InvRelabel}{\inf} d_{\square,\cF}(U^\psi,W) 
= \underset{\psi\in \Relabel}{\inf} d_{\square,\cF}(U^\psi,W) 	\\
& = \underset{\varphi, \psi\in \InvRelabel}{\inf} d_{\square,\cF}(U^\psi,W^\varphi) 
= \underset{\varphi,\psi\in \Relabel}{\min} d_{\square,\cF}(U^\psi,W^\varphi).
\end{aligned}
\end{equation}
\end{lemma}

\subsection{Equivalence of metrics and compactness}
In this section, we summarize two important results about probability graphons from \cite{abraham2023probabilitygraphons}.

The following theorem states the topological equivalence of the two unlabelled cut distances $\delta_{\square,\mathcal{LP}}$ and $\delta_{\square,\cF}.$ This is a slight reformulation of Corollary 5.6 and Remark 5.7 in \cite{abraham2023probabilitygraphons}.

%{\color{blue} \cite[Corollary~5.6]{abraham2023probabilitygraphons}}

\begin{theorem}[Corollary 5.6 and Remark 5.7 in \cite{abraham2023probabilitygraphons}]
	\label{cor:equiv-topo}
The cut distances $\delta_{\square,\mathcal{LP}}$ and $\delta_{\square,\cF},$ for every choice of the convergence determining sequence $\F$, induce the same topology  on  $\UGraphon$. 
\end{theorem}

The next theorem is a compactness result for families of probability graphons. This is a slight reformulation of Proposition 5.2 in \cite{abraham2023probabilitygraphons}.

%{\color{blue} \cite[Proposition~5.2]{abraham2023probabilitygraphons}}

\begin{theorem}[Proposition 5.2 in \cite{abraham2023probabilitygraphons}]\label{ThmRelCompact}
Let $\cK\subset \UKernel$ be a uniformly bounded subset, i.e. there exists $K>0$ such that for any $W\in \cK$$$\sup_{W\in \cK}\|W\|_{TV}<K.$$ Then the set $\cK$ is relatively compact for $\delta_{\square}$ (or equivalently $\delta_{\square,\cF}$) if and only if it is tight, i.e.\ the set of measures $\{ M_W : W\in\mathcal{K} \} \subset \Meas$ is tight, where $M_W$ for $W\in \UKernel$ is the measure
\begin{equation}
  \label{eq:def-MW}
M_W(\drv z) 
= \vert W \vert ([0,1]^2; \drv z)
= \int_{[0,1]^2} \vert W\vert (x,y;\drv z) \ \drv x \drv y.
\end{equation}
\end{theorem}
\begin{remark}
In particular, if $\Space$ is compact then $\UGraphon$ is compact. 
\end{remark}

\section{Equivalence Cut-metric and Quotients}\label{Sec4Results}

In this section, we generalise the notions of overlay functionals and quotient sets from real-valued graphons to the case of probability graphons, we study their properties and we show the main result of this work, Theorems \ref{ThmEquivalenceConvQuotientOverlay} and \ref{ThmEquivalenceConvQuotientOverlayVersion2}. The interested reader can find more details about overlay functionals and quotient sets for real-valued graphons in Chapter 12 in \cite{LovaszGraphLimits} or in \cite{borgs2011convergentAnnals}. Our proof of the main result follows the proof scheme of \cite{borgs2011convergentAnnals} and Chapter 12 in \cite{LovaszGraphLimits}.

\subsection{The overlay functional}
Let $\cC\subset \CbFunct$ and $\lambda$ the Lebesgue measure on $[0,1].$ A $\cC$ edge decorated graph (or simply $\cC-$graph) is a triple $G^{\beta}=(G,\beta(G))=(V(G),E(G),\beta(G))$ where $V=V(G)$ is the vertex set, $E(G)$ is the edge set, and $\beta(G)$ is a function
$$
\beta(G):V(G)\times V(G)\rightarrow \cC \cup \{0\}$$
such that $\beta(G)_{v,w}\neq 0$ if and only if $\{v,w\}\in E(G).$

For a probability measure $\alpha$ on $[k]$, we denote by $\Pi(\alpha)$ the set of partitions $\left(S_1, \ldots, S_k\right)$ of $[0,1]$ into $k$ measurable subsets with $\lambda\left(S_i\right)=\alpha_i$. For every $U\in \Kernel$ and $\CbFunct-$graph $G^{\beta}$ on the vertex set $[k]$, we define
\begin{equation}  \label{defiOverlay1}
\mathcal{C}(U, G^{\beta})=\sup _{\left(S_{ 1}, \ldots, S_k\right) \in \Pi(\alpha)} \sum_{i, j \in[k]}  \int_{S_i \times S_j} U(x, y,\beta_{i j}(G)) d x d y.
\end{equation}

Additionally to $\cC-$graphs, we define $\cC-$valued kernels for $\cC\subset \CbFunct $.

\begin{definition}[$\CbFunct$-valued kernels]\label{DefCBValuedKernel}
A \emph{$\cC$-valued kernel} is a map $W$ from $[0,1]^2$  to $\cC\subset\CbFunct$, 
such that:
\begin{enumerate}
\item   $W$  is  a continuous and bounded function belonging to $\cC$ in the $z\in \Space$ coordinate, i.e.\ for (almost) every $(x,y) \in [0,1]^2$, $W(x,y;\cdot)$ belongs to $\cC\subset \CbFunct$.
\item $W$ is \emph{stongly measurable} in $(x,y),$ i.e.\ there exists a sequence $(W_n)$ such that $W_n=\sum^{N_n}_if_{i,n}\mathbbm{1}_{A_{i,n}}$ where $A_{i,n}$ are measurable subsets of $[0,1]^2$ and $f_{i,n}\in \CbFunct $ (i.e.\ $W_n$ is a simple function) and $(W_n)$ converge to $W$ in $\|\cdot\|_{\infty}$ almost everywhere on $[0,1]^2$ with respect to the Lebesgue measure.
\item $W$ is \emph{bounded}, i.e.\
\begin{equation}
  \label{eq:def:TMCb}
\sup\|W\|_{\infty}:=    \sup_{x,y\in [0, 1] }\, \TotalMass{W(x, y; \cdot)} <+\infty 
  \end{equation}where the supremum is meant to be an essential supremum (almost everywhere with respect to the Lebesgue measure).

\end{enumerate}
\end{definition}

Note that in Definition \ref{DefCBValuedKernel} we assumed strong measurability (sometimes also called Bochner measurability) which is stronger than (Borel) measurability in general. Pettis measurability theorem states that a function is strongly measurable if and only if it is measurable and its image is separable, see for example Page 16 in \cite{StrongMeasEquivMeasSeparNBanach}. In particular, when $\CbFunct$ is separable, i.e.\ when $\Space$ is compact, the two notions are equivalent.

In the following, for a $\CbFunct$-valued kernel $W$ and a measure-preserving map $\varphi\in \Relabel,$ we will denote by $W^{\varphi}$ the $\CbFunct$-valued kernel $W^{\varphi}(x,y)=W(\varphi(x),\varphi(y)).$ Similarly to the case of probability graphons we will say that two $\CbFunct$-valued kernels $W$ and $U$ are \emph{weakly isomorphic} if $W^{\varphi}=U^{\psi}$ almost everywhere for some $\varphi,\psi \in \Relabel.$ Moreover, we will denote by $\widetilde{\cU}$ the space of equivalence classes of $\CbFunct$-valued kernels up to weak isomorphism.
    
Similarly to signed measure-valued kernels, we define step-functions for $\CbFunct-$kernel.

\begin{definition}[$\CbFunct$-valued stepfunctions]
  \label{def:stepfunctionCb}
  A $\CbFunct$-valued kernel $W$ is a \emph{step-function} if there exists a
  finite partition of $[0,1]$ into measurable (possibly empty) sets, say
  $\mathcal{P}=\{S_1,\cdots,S_k\}$, such that $W$ is constant on the sets
  $S_i \times  S_j$, for $1\leq i,j\leq  k$.  
\end{definition}

The following lemma asserts that $\CbFunct$-valued step-functions are enough to approximate $\CbFunct$-valued kernels.

\begin{lemma}\label{LemmApproxCbkernelsGeneral}
Let $U$ be a $\CbFunct-$valued kernel. For every $\varepsilon>0$ there exists a $\CbFunct-$valued step-function $U^{\prime}$ such that 

$$
\int_{[0,1]^2}\|U-U^{\prime}\|_{\infty}<\varepsilon.
$$
\end{lemma}
 \proof 
By the Bochner measurability assumption in the definition of $\CbFunct-$valued kernels (Definition \ref{DefCBValuedKernel}), there exists a simple $\CbFunct-$function arbitrarily close to $W$ almost everywhere. Therefore, it suffices to show that we can approximate a simple $\CbFunct-$function with a $\CbFunct-$valued step-function almost everywhere. The proof follows now simply by recalling that the Lebesgue measure on $[0,1]^2$ is the product measure of the Lebesgue measure on $[0,1]$ and standard measure-theoretic arguments. 
 \endproof

 We now explain how $\cC-$graphs and $\cC-$valued kernels are related by  $\cC-$valued step-functions.
A $\cC-$graph can be naturally represented as a $\cC-$valued step-function in the following way:

Let $G^{\beta}=(G,\beta(G))=(V(G),E(G),\beta(G))$ be a $\CbFunct-$graph with $|V(G)|=n.$
We will label the vertices $V(G)=V=\{v_1,\ldots,v_n\}.$ 
We will consider the interval $[0,1]$ and partition it into $n$ intervals of length $1/n$. For $i\in [n]$ we let $J_i=((i-1)/n,i/n].$ We then define $W_G$ by:

\begin{equation}\label{eq:DefWG}
    W_{G^{\beta}}(x,y)=\begin{cases}
 \beta_{i,j}(G)=\beta_{v_i,v_j}(G)\text{ if }(x,y)\in J_i\times J_j\subset [0,1]^2 \text{ and } \{v_i,v_j\}\in E(G)
 \\ 

0 \text{ else }
   
\end{cases}
\end{equation}
We state now another approximation result, that can be interpreted as a refinement of Lemma \ref{LemmApproxCbkernelsGeneral} for the approximation of $\CbFunct-$kernels taking values in a finite subspace of $\CbFunct.$ 

\begin{lemma}\label{lemmaApproxDecGraphonsRealValGraphons}
Let $U(x,y;z)=\sum^N_{i=1}g_i(z) u_i(x,y)$ be a $\CbFunct-$valued kernel, where $u_i$ are real-valued kernels, such that $u_i\in L^1([0,1]^2)$ and $g_i\in \CbFunct$. Then for every $\varepsilon>0$ there exist $w_i$ step real-valued kernels (composed of rectangles) such that  
$$
   \int_{[0,1]^2}\left\|U- W\right\|_{\infty}\leq \varepsilon $$
where $W(x,y)=\sum^N_{i=1}g_i w_i(x,y).$
\end{lemma}
\proof
Let's choose $w_i$ step real-valued kernels (composed of rectangles) such that  
$$
\|u_i-w_i\|_{1} \leq \frac{\varepsilon}{N\max_{i\in [N]}\|g_i\|_{\infty}}.
$$
 We know that such real-valued step-functions exist for real-valued kernels in $L^1$. Therefore, we obtain the result
$$
\begin{aligned}
  & \int_{[0,1]^2}\left\|U- W\right\|_{\infty}=\int_{[0,1]^2}\left\|\sum^N_{i=1}g_i u_i(x,y)- \sum^N_{i=1}g_i w_i(x,y)\right\|_{\infty} \leq  \\
  &\sum^N_{i=1}\|g_i\|_{\infty}\int_{[0,1]^2}|u_i(x,y)-w_i(x,y)|\leq N \max_{i\in [N]}\|g_i\|_{\infty}\frac{\varepsilon}{N\max_{i\in [N]}\|g_i\|_{\infty}}=\varepsilon,
   \end{aligned}$$
where in the first inequality we used the absolute homogeneity of norms.
\endproof

To define weak convergence of probability measures one considers the application of a measure $\mu$ to a function $f\in\CbFunct,$ obtained integrating $f$ with respect to $\mu,$ that we denoted with $\mu[f].$ Similarly, one can naturally define the application of $U\in \UKernel,$ a signed measure-valued kernel, to $W$ a $\CbFunct-$valued kernel. We denote the application of $U$ to $W$ with $U(W),$ the real-valued kernel such that $$U(W)(x,y)=U(x,y)(W(x,y))=U(x,y;W(x,y;\cdot))=\int_{\Space}W(x,y;z)U(x,y,\drv z).$$

We show now that this operation is well defined.

\begin{theorem}
For $U\in \UKernel$ and $W$ a $\CbFunct-$valued kernel the function 

$$U(W)(x,y)$$
from $[0,1]\times [0,1]$ to $\R$ is measurable, i.e. $U(W)$ is a real-valued kernel.
\end{theorem}
\proof
Let's assume first that $W$ is a measurable step-function $W=\sum^N_{i,j=1}f_{ij} \mathbbm{1}_{P_i\times P_j}$ where $f_{ij} \in \CbFunct $  and $P_1,\ldots,P_N$ is a measurable partition of $[0,1]^2.$ Thus, we have,

$$
\begin{aligned}&U(W)(x,y)=U(x,y)(W(x,y))\\
&=U(x,y;W(x,y;\cdot))=\int_{\Space}W(x,y;z)U(x,y,\drv z)\\
&= \int_{\Space}\sum^N_{i,j=1}f_{ij}(z) \mathbbm{1}_{P_i\times P_j}(x,y)U(x,y,\drv z)\\
&=\sum^N_{i,j=1} \mathbbm{1}_{P_i\times P_j}(x,y)\int_{\Space}f_{ij}(z) U(x,y,\drv z)\\
&=\sum^N_{i,j=1} \mathbbm{1}_{P_i\times P_j}(x,y) U(x,y,f_{i,j})
\end{aligned}$$
Therefore, $U(W)$ is measurable as $U(x,y,f_{i,j})$ is measurable by the measurability condition in the definition of probability graphon (Definition \ref{Def:MeasValuedKernel}) and multiplication and sum of measurable (real-valued) functions is measurable.

Let's now consider a general Bochner measurable $\CbFunct-$kernel. By the definition of Bochner measurability and Lemma \ref{LemmApproxCbkernelsGeneral}, there exists a sequence $W_n$ of measurable $\CbFunct-$valued step-functions converging pointwise to $W$, i.e.\ $\|W_n(x,y,\cdot)-W(x,y,\cdot)\|_{\infty}=\sup_{z\in \Space}|W_n(x,y,z)-W(x,y,z)|$ tends to $0$ for every $(x,y)\in [0,1]^2\setminus N$ where $N\subset[0,1]^2$ has Lebesgue measure $0.$  By the dominated convergence theorem, we obtain
$$
\begin{aligned}
&\lim_{n\rightarrow \infty}U(W_n)(x,y)\\
&=\lim_{n\rightarrow \infty}\int_{\Space}W_n(x,y;z)U(x,y,\drv z)\\
&=\int_{\Space}W(x,y;z)U(x,y,\drv z)\\
&= U(W)(x,y).
\end{aligned}
$$
We could apply the dominated convergence theorem because for every $(x,y)\in[0,1]^2$ we have that $|W_n(x,y,\cdot)|$ is asymptotically dominated by $\sup\|W\|_{\infty}+1$ that is integrable with respect to $U(x,y,\mathrm{d}z)$ as $U$ is a probability graphon. 

Therefore, it follows that $U(W)$ is measurable as the pointwise limit of measurable functions. 
\endproof

We now proof an identity for $U(W)$ in the special case of $\CbFunct-$valued kernel $W$ taking values in a finite dimensional subspace of $\CbFunct.$

\begin{lemma}\label{ProprLinearProbGraphons}
Let $w_1,\ldots, w_N$ be real-valued graphons $f_1,\ldots , f_N\in\CbFunct$ and $U$ be a probability graphon.  Then $U(\sum^{N}_{n=1}f_n w_n)=(\sum^{N}_{n=1} w_nU[f_n])$ almost everywhere on $[0,1]^2.$
\end{lemma}
\proof
Let $w_1,\ldots, w_N$ be real-valued graphons and $f_1,\ldots,f_N \in \CbFunct$ continuous and bounded functions. The result follows directly from the simple almost everywhere equality on $[0,1]^2$
$$\begin{aligned}
U(\sum^{N}_{n=1}f_n w_n)(x,y)&= U(x,y)(\sum^{N}_{n=1}f_n w_n(x,y))\\
&=\int_{\Space}\sum^{N}_{n=1}f_n(z) w_n(x,y)U(x,y,\mathrm{d}z)
\\
& =\sum^{N}_{n=1}w_n(x,y)\int_{\Space}f_n(z) U(x,y,\mathrm{d}z)
\\ &=\sum^{N}_{n=1} w_n(x,y)U[f_n](x,y)\\ &=(\sum^{N}_{n=1} w_nU[f_n])(x,y),
\end{aligned}
$$
where we used the linearity of the integral in the third equality.
\endproof

We can extend the previous theorem to infinite sums when the monotone convergence theorem or the dominated convergence theorem holds as we show in the following lemma.

\begin{lemma}
    Let $(w_n)_n$ be a sequence of non-negative real-valued kernels, $(f_n)_n$ a sequence in $\CbFunct$ and $U$ a probability graphon. Then $U(\sum^{\infty}_{n=1}f_n w_n)=(\sum^{\infty}_{n=1} w_nU[f_n])$ almost everywhere on $[0,1]^2.$
\end{lemma}
\proof
By monotone convergence (the function $f_nw_n$ is non-negative as $w_n$ is non-negative and $f_n$ is $[0,1]-$valued as an element of a convergence determining sequence) we have that it exists a set $N\subset [0,1]^2$ of measure zero such that for every $(x,y)\in [0,1]^2\setminus N$
$$
\begin{aligned}
&\sum^{\infty}_{n=1}w_n(x,y)U[f_n](x,y) =\sum^{\infty}_{n=1}w_n(x,y) \int_{\Space}f_n(z)U(x,y,\mathrm{d}z)\mathrm{d} x \mathrm{d} y\\
&= \int_{\Space}\sum^{\infty}_{n=1}  f_n(z)w_n(x,y)U(x,y,\mathrm{d}z)= U(\sum^{\infty}_{n=1} f_n w_n)(x,y).
\end{aligned}
$$
This concludes the proof.
\endproof

We can now further generalize the definition of overlay functional \eqref{defiOverlay1} and define, for $U\in \Kernel$ and $W$ a $\CbFunct-$valued kernel $W,$
\begin{equation}\label{Eq:DefOverFunc}
\mathcal{C}(U, W)=\sup _{\varphi \in S_{[0,1]}} \int_{[0,1]^2}U(W^{\varphi})(x,y) \mathrm{d}x \mathrm{d}y=\sup _{\varphi \in S_{[0,1]}} \int_{[0,1]^2} U(x, y) (W(\varphi(x), \varphi(y))) \mathrm{d}x \mathrm{d}y. 
\end{equation}
\begin{remark}
There are no measurability issues. Observe in fact that if $W$ is Bochner measurable $W^{\varphi}(x,y)=W(\varphi(x), \varphi(y))$ is still Bochner measurable. This follows simply from the fact that if $W$ is a step-function then $W^{\varphi}$ is also a step-function.
\end{remark}

Equation \eqref{Eq:DefOverFunc} extends the definition \eqref{defiOverlay1} in the sense that if $U$ is any probability graphon and $G$ is a $\CbFunct-$graph, then
\begin{equation}\label{eq:IdentOverlayGraphGraphon}
\mathcal{C}(U, G^{\beta})=\mathcal{C}\left(U, W_{G^{\beta}}\right),
\end{equation}
where $W_{G^{\beta}}$ is the $\CbFunct-$kernel obtained by $\CbFunct-$graph $G^{\beta}$, recall \eqref{eq:DefWG}.

We will call  $\mathcal{C}(U, W)$ the \emph{overlay functional}. We will now state a few properties of this functional. Similarly to the statement for distances Lemma \ref{thm_min_dist} (Proposition 3.18 in \cite{abraham2023probabilitygraphons}), we have that
$$
\begin{aligned}
\mathcal{C}(U, W) & =\sup _{\varphi \in S_{[0,1]}} \int_{[0,1]^2}U( W^{\varphi})=\sup _{\varphi \in S_{[0,1]}}\int_{[0,1]^2} U^{\varphi}( W)=\sup _{\varphi, \psi \in \bar{S}_{[0,1]}} \int_{[0,1]^2}U^{\varphi}( W^\psi) \\
& =\sup \left\{ \int_{[0,1]^2}U_0(W_0) :\left(\exists \varphi, \psi \in \bar{S}_{[0,1]}\right) U=U_0^{\varphi}, W=W_0^\psi\right\} .
\end{aligned}
$$

Therefore, the overlay functional is invariant under measure preserving transformations of the finite signed measure-valued kernels and the $\CbFunct-$kernels, i.e., it is also a functional on the space $\UKernel \times \widetilde{\cU}$ (or its subsets) where we recall that $\UKernel$ is the set of isomorphism classes of weakly isomorphic signed measure-valued kernels and $\widetilde{\cU}$ is the set of equivalence classes of weakly isomorphic $\CbFunct-$kernels. 

Moreover, the functional $\mathcal{C}(U, W)$ is sub-additive in each variable, i.e.\ for $U,V \in \UKernel$ and $W,Q$ $\CbFunct-$valued kernels, we have
\begin{equation}\label{subaddOverlay}
    \mathcal{C}(U+V, W) \leq \mathcal{C}(U, W)+\mathcal{C}(V, W)
\end{equation}

and 

\begin{equation}\label{subaddOverlay2}
    \mathcal{C}(U, W+Q) \leq \mathcal{C}(U, W)+\mathcal{C}(U, Q).
\end{equation}

We remark that the overlay functional is only sub-additive and not bi-linear.

However, for $\lambda\in \R$ we have $$
\mathcal{C}(\lambda U,W)=\mathcal{C}(U, \lambda W)
$$
and\ if $\lambda>0$, then
\begin{equation}\label{Poslinearityoverlay}
\mathcal{C}(\lambda U,W)=\mathcal{C}(U, \lambda W)=\lambda \mathcal{C}(U, W),   
\end{equation}
i.e.\ the overlay functional is homogeneous for positive scalars. Moreover, we have  $\mathcal{C}(U, W)=\mathcal{C}(-U,-W)$. However, we remark that, in general,  $\mathcal{C}(U, W)$ and $\mathcal{C}(-U, W)$ are unrelated.\newline

We show now a result from functional analysis that we will need in the following.\newline
Let $\cX$ be a normed vector space (over $\R$). For a (possibly infinite) set $\cC\subset \cX$ we will denote by $span(\cC)$ the span of $\cC$, i.e. the vector space given by the finite linear combinations of elements of $\cC.$ In other words, 
$$
\begin{aligned}
span(\cC)=&\{g\in  \CbFunct \text{ such that there exists }N\in\N, \ a_1\ldots a_N\in \R \text{ and }f_1,\ldots,f_N\in \cC \\
&\text{ such that } y=\sum^N_{i=1}a_if_i\}.
\end{aligned}$$

\begin{lemma}\label{lemmaSeparabSpanCount}
Let $\cX$ be a normed vector space (over $\R$).  For a countable sequence $(f_i)_{i\in \N}=\cC\subset \cX,$ the closure of $span(\cC)$, denoted with $\overline{span(\cC)},$ is separable in $\cX$. 
\end{lemma}
\proof
We observe that, for every $N\in\N,$ the set 
$$
\begin{aligned}
\cB_N=&\{g\in  \CbFunct \text{ such that there exists } a_1\ldots a_n\in \Q \text{ such that } y=\sum^N_{i=1}a_if_i\}.
\end{aligned}$$
is countable.

Therefore, the set \begin{equation}\label{eq:CountabSet}\begin{aligned}
&\cB=\cup_{N\in \N}\cB_N=\\
&\{g\in  \CbFunct \text{ such that there exists } N\in\N \text{ and }a_1\ldots a_n\in \Q \text{ such that } y=\sum^N_{i=1}a_if_i\}\end{aligned}\end{equation} is also a countable set as a countable union of countable sets is countable.
Moreover, we have that $\cB$ is dense in $\overline{span(\cC)}$ and thus $\overline{span(\cC)}$ is separable.
\endproof

We show now that the overlay functional is continuous in the first variable with respect to the metric $\delta_{\square,\cF}.$

In the following lemma, $\cF$ is a countable sequence in $\CbFunct$. Therefore, with $\overline{span(\cF)}$ we mean the closure of the span of $\cF$ with respect to the topology of $\CbFunct$. This is a closed subspace of $\CbFunct.$

\begin{lemma}\label{LemmmetricConvImpliesOverlay}
Let $(f_{n})_{n\in \N}=\cF$ be a convergence-determining sequence. Let $U \in \Kernel$ and $(U_n)_n$ be a sequence such that $U_n\in \Kernel.$ If $\delta_{\square,\cF}\left(U_n, U\right) \rightarrow 0$ as $n \rightarrow \infty$, then for every $\overline{span(\cF)}$-valued kernel $W$ we have $\mathcal{C}\left(U_n, W\right) \rightarrow \mathcal{C}(U, W)$.
\end{lemma}
\proof
By subadditivity \eqref{subaddOverlay}, we have
$$
 -\mathcal{C}\left(U-U_n, W\right) \leq \mathcal{C}\left(U_n, W\right)-\mathcal{C}(U, W) \leq \mathcal{C}\left(U_n-U, W\right).
$$
Thus we only have to show that $\mathcal{C}\left(U_n-U, W\right), \mathcal{C}\left(U-U_n, W\right) \rightarrow 0$.
As $\cC(U,W)$ is invariant up to weak isomorphism of $U$ we can assume that $U_n$ are optimally overlayed with $U$. Therefore, $\delta_{\square,\cF}\left(U_n-U,0\right)=\|U_n-U\|_{\square,\cF}=\delta_{\square,\cF}\left(U_n, U\right)\rightarrow 0.$
In other words, it is enough to show the lemma in the case when $U=0$.

By definition of overlay functional, we have $\mathcal{C}(U_n,W)\geq \int_{[0,1]^{2}}U_n(W)\geq -\left|\int_{[0,1]^{2}}U_n(W)\right|,$ and we will prove now that the right-hand side converges to 0. 

Let's show first that $$
\left|\int_{[0,1]^{2}}U_n(W)\right|\rightarrow 0
$$
when $W$ is an indicator function of a rectangle, i.e.\ $W=f\mathbbm{1}_{S\times T}$ for a suitable class of functions $f\in \overline{span(\cF)}$ and measurable sets $S,T\subset [0,1].$ 

Let's define
\begin{equation*}\begin{aligned}
&\cB=\{g\in  \CbFunct \text{ such that there exists } N\in\N \text{ and }a_1\ldots a_n\in \Q \text{ such that } y=\sum^N_{i=1}a_if_i\}.\end{aligned}\end{equation*}
From Lemma \ref{lemmaSeparabSpanCount} and its proof (see \eqref{eq:CountabSet}) the set $\cB$ is a countable subset of  $span(\cF)$ that is dense in $\overline{span(\cF)}.$ %{\color{blue} Maybe $\cF$ is enough.......}
Therefore, let $f\in \cB$, i.e.\ $f=\sum^N_{i=1}a_if_{m_i}$ for $a_i\in \Q$ and $f_{m_i}\in (f_k)_{k\in \N}=\cF$ and $N\in \N$, then
$$\left|\int_{[0,1]^{2}}U_n(W)\right|=\left|\int_{[0,1]^{2}}U_n[f]\mathbbm{1}_{S\times T}\right|\leq \|U_n[f]\|_{\square,\R}\leq\sum^N_{i=1}|a_i| \|U_n[f_{m_i}]\|_{\square,\R}\leq a 2^m \|U_n\|_{\square,\cF}\rightarrow 0.$$
where $a=\sum^N_{i=1}|a_i|$ and $m=\max_{i\in [N]}\{m_i\}$ and we used Lemma \ref{lemmaBoundConvDetSequenceCutNormF} in the last inequality.

By linearity, the same holds if $W$ is a $\cB-$valued step-function.

Let's suppose that $W$ is now a general $\overline{span(\cF)}$-valued kernel. We can approximate $W$ with a $\cB-$valued step-function $V,$ i.e.\ for any $\varepsilon>0$ there exists a $\cB-$valued step function $V$ such that $\int_{[0,1]^2}\|W-V\|_{\infty}\leq \varepsilon.$ This is possible as $\cB$ is dense in $\overline{span(\cF)}$ and $V$ is Bochner measurable (that is equivalent to measurable in this case by Pettis measurability theorem as $\overline{span(\cF)}$ is separable).  We obtain the following bounds:
$$
\begin{aligned}
   \left|\int_{[0,1]^{2}}U_n(W)\right|\leq 
   &\left|\int_{[0,1]^{2}}U_n(W-V)\right|+\left|\int_{[0,1]^{2}}U_n(V)\right|\leq \\
   &\sup\|U_n\|_{TV}\int_{[0,1]^2}\|W-V\|_{\infty}+\left|\int_{[0,1]^{2}}U_n(V)\right|\leq \\
     & \varepsilon \sup\|U_n\|_{TV}+\left|\int_{[0,1]^{2}}U_n(V)\right|.
\end{aligned}$$
Therefore, $\int_{[0,1]^{2}}U_n(W)\rightarrow 0$ as $\varepsilon$ can be chosen arbitrarily small and $|\int_{[0,1]^{2}}U_n(V)|\rightarrow 0$ as $V$ is a step function. Hence $\liminf_n\mathcal{C}(U_n,W)=0.$ \newline

To prove the opposite inequality, we start with the case when $W$ is a $span(\cF)-$valued step-function, i.e.\ $W=\sum^m_{i=1}g_i\mathbbm{1}_{S_i\times T_i}$ where $S_i,T_i$ are measurable subsets of $[0,1]$ and $g_i=\sum^{N_i}_{j=1}a_{ij}f_{m_{i,j}}$ for $a_{ij}\in \R$ and $f_{m_{i,j}}\in \cF.$ We obtain
$$
\begin{aligned}
\mathcal{C}\left(U_n, W\right) & \leq \sum_{i=1}^m \mathcal{C}\left(U_n, g_i \mathbbm{1}_{S_i \times T_i}\right)=\sum_{i=1}^m\sup_{\phi}\int_{[0,1]^2}U_n(g_i \mathbbm{1}^{\phi}_{S_i \times T_i}) = \\ & \sum_{i=1}^m\sup_{\phi}\int_{[0,1]^2}U_n[g_i] \mathbbm{1}^{\phi}_{S_i \times T_i}\leq\\
&  \sum_{i=1}^m\left\|U_n[g_i]\right\|_{\square,\R}\leq \sum_{i=1}^m\sum_{j=1}^{N_i}|a_{i,j}|\left\|U_n[f_{m_{i,j}}]\right\|_{\square,\R}\leq a2^m \left\|U_n\right\|_{\square,\cF}\rightarrow 0 
\end{aligned}
$$
where $a=\sum_{i=1}^m\sum_{j=1}^{N_i}|a_{i,j}|$ and $m=\max_{i\in [m],j\in [N_i]}\{m_{i,j}\}$ and we used Lemma \ref{lemmaBoundConvDetSequenceCutNormF} in the last inequality.

%The $ \|\cdot\|_{\infty} $ is the total variation for measures and infinity norm for continuous and bounded functions, also below.

%Now if $W=\sum^N_{i=1}g_iw_i(x,y)$ for arbitrary real-valued $L^1$ kernels, then for every $\varepsilon>0$ we can find a .......
Let's consider now $W$ to be a general $\overline{span(\cF)}-$valued kernel. Then for every $\varepsilon>0$ we can find a $span(\cF)-$valued step-function $W^{\prime}$ such that $\sup\left\|U_n\right\|_{TV}\int_{[0,1]^2}\left\|W-W^{\prime}\right\|_{\infty}  \leq \varepsilon / 2.$

We know that $\mathcal{C}\left(U_n, W^{\prime}\right) \rightarrow 0$, and hence $\mathcal{C}\left(U_n, W^{\prime}\right) \leq$ $\varepsilon / 2$ if $n$ is large enough. But then, for $n$ large enough, 
$$
\mathcal{C}\left(U_n, W\right) \leq \mathcal{C}\left(U_n, W-W^{\prime}\right)+\mathcal{C}\left(U_n, W^{\prime}\right) \leq\sup\left\|U_n\right\|_{\infty}\int_{[0,1]^2}\left\|W-W^{\prime}\right\|_{\infty} +\varepsilon / 2 \leq \varepsilon .
$$

Therefore, we obtained that $\lim \sup _n \mathcal{C}\left(U_n, W\right) \leq 0$, and this completes the proof of the lemma.
\endproof

A direct consequence of the previous lemma is the following corollary.

\begin{corollary}\label{CorSequencesEqOv}
 Let $U \in \UKernel$ and $(U_n)_n$ be a sequence such that $U_n\in \UKernel.$ If $\delta_{\square,\mathcal{LP}}\left(U_n, U\right) \rightarrow 0$ as $n \rightarrow \infty$, then for every (Bochner measurable) $\CbFunct$-valued kernel $W$ we have $\mathcal{C}\left(U_n, W\right) \rightarrow \mathcal{C}(U, W)$.
\end{corollary}
\proof
This follows from Lemma \ref{LemmmetricConvImpliesOverlay}, Theorem \ref{cor:equiv-topo} (Corollary 5.6 in \cite{abraham2023probabilitygraphons}) and the observation that $\delta_{\square,\mathcal{LP}}$ is independent by a specific convergence determining sequence $\cF.$
\endproof

We also prove the continuity of the overlay functional with respect to particular sequences of $\CbFunct-$kernels.

\begin{lemma}
   \label{continuityContBoundGraphRealValCutDist}
   Let $U=\sum^N_{i=1}f_iu_{i}$ and $(U_n)_n$ be a sequence such that $U_n=\sum^N_{i=1}f_iu_{i,n}$ for $f_i\in \CbFunct$ and $u_1,\ldots,u_N,u_{1,n},\ldots,u_{N,n}$ real-valued graphons. If \begin{equation}\label{Eq:CutmetricFiniteCb}
    \delta_{\square,N}(U,U_n)=\inf_{\phi}\sum^N_{i=1}\NcutR{u_i-u^\phi_{i,n}}\rightarrow0
\end{equation} as $n \rightarrow \infty$, then for every $W \in\UKernel$ we have $\mathcal{C}\left( W,U_n\right) \rightarrow \mathcal{C}(W,U)$.
\end{lemma}

\proof
Let $U=\sum^N_{i=1}f_iu_{i}$ and the sequence $U_n=\sum^N_{i=1}f_iu_{i,n}$ as in the statement of the lemma.
By subadditivity \eqref{subaddOverlay}, we have
$$
 -\mathcal{C}\left(W,U-U_n \right) \leq \mathcal{C}\left(W,U_n\right)-\mathcal{C}(W,U) \leq \mathcal{C}\left(W,U_n-U\right).
$$
Hence it is enough to show that $\mathcal{C}\left(W,U_n-U\right)$, $\mathcal{C}\left(W,U-U_n\right) \rightarrow 0$. This means we need to prove the statement only when $U=0,$ i.e.\ when $u_1,\ldots,u_N=0.$

By definition, we have $\mathcal{C}(W,U_n)\geq \sum^N_{i=1}\int_{[0,1]^{2}}W[f_i]u_{i,n}\geq-\sum^N_{i=1}\left|\int_{[0,1]^{2}}W[f_i]u_{i,n} \right|$, and the right-hand-side tends to $0$ by Lemma 8.22 in \cite{LovaszGraphLimits}. Hence $\liminf_n\mathcal{C}(U_n,W)\geq 0.$\newline

We now show the other inequality to conclude the proof. We assume first that $W$ is a step-function, i.e.\ $W=\sum_{i=1}^m \mu_i \mathbbm{1}_{P_i \times Q_i}$ for $\mu_i$ signed measures and $P_i,Q_i$ measurable subsets of $[0,1].$ In this case, the result follows from the following inequalities: 

$$
\begin{aligned}
\mathcal{C}\left(W,U_n\right) & \leq \sum_{i=1}^m \mathcal{C}\left( \mu_i \mathbbm{1}_{P_i \times Q_i},U_n\right)=\sum_{i=1}^m\sup_{\phi}\int_{[0,1]^2}\mu_i \mathbbm{1}_{P_i \times Q_i}(U^{\phi}_n) = \\ & \sum_{i=1}^m\sup_{\phi}\int_{[0,1]^2}\sum^N_{j=1}\mu_i(f_j ) u^{\phi}_{j,n}\mathbbm{1}_{P_i \times Q_i}=\\
& 
\sum_{i=1}^m\sum^N_{j=1}\mu_i(f_j )\sup_{\phi}\int_{P_i \times Q_i} u^{\phi}_{j,n} \leq \\
&
\sum_{i=1}^m\sum^N_{j=1}|\mu_i(f_j )|\|u_n\|_{\square,\R}
\rightarrow 0 .
\end{aligned}
$$
Let's now consider the general case $W\in \UKernel.$ From Lemma 4.5 in \cite{abraham2023probabilitygraphons} we know that for every $\varepsilon>0$ it exists a step-valued kernel $W^{\prime}\in \UKernel$ such that $\sum^N_{i=1}\|u_{i,n}\|_{\infty}\int_{[0,1]^2}|W[f_i]-W^{\prime}[f_i]|<\varepsilon/2.$
Moreover, we know that $\mathcal{C}\left(W^{\prime},U_n\right) \rightarrow 0$, and hence $\mathcal{C}\left(W^{\prime},U_n\right) \leq$ $\varepsilon / 2$ if $n$ is large enough. But then

$$\begin{aligned}
\mathcal{C}\left( W,U_n\right) \leq \mathcal{C}\left(W-W^{\prime},U_n\right)+\mathcal{C}\left( W^{\prime},U_n\right) \leq\sup_{\varphi}\int_{[0,1]^2}(W-W^{\prime})(U^{\varphi}_n)+\varepsilon / 2 \leq \varepsilon ,
\end{aligned}
$$
where the last inequality follows from the inequality:
$$\begin{aligned}
&\sup_{\varphi}\int_{[0,1]^2}(W-W^{\prime})(U^{\varphi}_n)=\sup_{\varphi}\int_{[0,1]^2}(W-W^{\prime})(\sum^N_{i=1}f_iu^{\varphi}_{i,n})=\\
&\sup_{\varphi}\sum^N_{i=1}\int_{[0,1]^2}(W[f_i]-W^{\prime}[f_i])u^{\varphi}_{i,n}\leq \sum^N_{i=1}\|u_{i,n}\|_{\infty}\int_{[0,1]^2}|W[f_i]-W^{\prime}[f_i]|\leq \varepsilon/2.
\end{aligned}
$$
This shows that $\lim \sup _n \mathcal{C}\left( W,U_n\right) \leq 0$, and completes the proof.
\endproof

We have that the overlay functional $\cC$ is continuous in the second variable for general $\CbFunct-$valued kernels with respect to other norms. The following lemma is an example.

\begin{lemma}\label{continuityContBoundGraphL1}
Let $U$ be a $\CbFunct-$valued kernel and , $(U_n)$ a sequence of $\CbFunct-$valued kernels that are Bochner measurable. If $\int_{[0,1]^2}\left\|U_n- U\right\|_{\infty} \rightarrow 0$ as $n \rightarrow \infty$ then for every $W $ measure valued graphon we have $\mathcal{C}\left( W,U_n\right) \rightarrow \mathcal{C}(W,U)$.
\end{lemma}
\proof
By subadditivity \eqref{subaddOverlay2}, we have
$$
-\mathcal{C}\left(W,U-U_n\right) \leq \mathcal{C}\left(W,U_n\right)-\mathcal{C}(W,U) \leq \mathcal{C}\left(W,U_n-U\right),
$$
and hence it is enough to prove that $\mathcal{C}\left(W,U_n-U\right), \mathcal{C}\left(W,U-U_n, \right) \rightarrow 0$. In other words, we only need to show the lemma in the case when $U=0$.

By definition, we have $$\sup\|W\|_{\infty}\int_{[0,1]^{2}}\|U_n\|_{\infty}\geq \mathcal{C}(W,U_n)\geq -\sup\|W\|_{\infty}\int_{[0,1]^{2}}\|U_n\|_{\infty}.$$ By assumption, we have $$
\int_{[0,1]^{2}}\|U_n\|_{\infty}\rightarrow 0
$$
and, therefore, we obtain that $\mathcal{C}(U_n,W)$ tends to $0.$
\endproof

\subsection{Another overlay functional}
In the following, we will need another quantity related to the overlay functional.  Let $\mathcal{F}=(f_n)_n$ be a convergence-determining sequence. We define the $\cF-$\emph{overlay functional} for two signed measure-valued kernels $U,W\in\UKernel$ as
$$
\mathcal{C}_{\mathcal{F}}(U,W)=\sup_{{\varphi \in S_{[0,1]}}}\sum^{\infty}_{n=0}2^{-n}\langle U[f_n],W^{\varphi}[f_n]\rangle=\sup_{{\varphi \in S_{[0,1]}}}\sum^{\infty}_{n=0}2^{-n}\int_{[0,1]^2}U[f_n]W^{\varphi}[f_n].
$$
This suggests that $\mathcal{C}_{\mathcal{F}}(.,.)$ can be interpreted in some sense as an inner product. To make this analogy more precise let's define the metric $\delta_{2,\mathcal{F}}$ derived from the $L^2$-norm as follows: 

$$
\delta_{2,\mathcal{F}}(U,W)=\inf_{ {\varphi \in S_{[0,1]}}}\|W-U^{\phi}\|_{2,\mathcal{F}}=\inf_{ {\varphi \in S_{[0,1]}}}\sum^{\infty}_{n=0}2^{-n}\|W[f_n]-U^{\phi}[f_n]\|_2,
$$
where $U,W\in \UKernel$ are signed measure valued kernels. 

We observe the following identity, reminiscent of the cosine theorem for scalar products:
\begin{equation}\label{EqOverlayFunctConvDetFam}
\begin{aligned}
\mathcal{C}_{\mathcal{F}}(U, W) & =\frac{1}{2}\left(\|U\|_{2,\mathcal{F}}^2+\|W\|_{2,\mathcal{F}}^2-\delta_{2,\mathcal{F}}(U, W)^2\right) \\
& =\frac{1}{2}\left(\delta_{2,\mathcal{F}}(U, 0)^2+\delta_{2,\mathcal{F}}(W, 0)^2-\delta_{2,\mathcal{F}}(U, W)^2\right) .
\end{aligned}
\end{equation}
This follows from the chain of equalities:
$$
\begin{aligned}
\delta_{2,\cF}(U, W)^2 & =\inf _{\varphi \in S_{[0,1]}}\left\|U-W^{\varphi}\right\|_{2,\mathcal{F}}^2=\|U\|_{2,\mathcal{F}}^2+\|W\|_{2,\mathcal{F}}^2-2 \sup _{\varphi \in S_{[0,1]}}\left\langle U, W^{\varphi}\right\rangle_{\mathcal{F}} \\
& =\|U\|_{2,\mathcal{F}}^2+\|W\|_{2,\mathcal{F}}^2-2 \mathcal{C}_{\mathcal{F}}(U, W) .
\end{aligned}
$$ 

We want to remark that the $\cF-$overlay functional $\cC_{\cF}$ takes as inputs two finite signed measure-valued kernels, differently from the overlay functional $\cC$ which takes one signed measure valued kernel and a $\CbFunct-$ valued kernel. 

However, the $\cF-$overlay functional shares a lot of properties with the overlay functional as sub-linearity and (positive) homogeneity.  
Moreover, also for the $\cF-$overlay functional we have the chain of equalities:
\begin{equation}\label{EqWeakIsoFOverlaFunc}
    \begin{aligned}
\mathcal{C}_{\cF}(U, W) & =\sup _{\varphi \in S_{[0,1]}} \int_{[0,1]^2}U[f_n] W^{\varphi}[f_n]=\sup _{\varphi \in S_{[0,1]}}\int_{[0,1]^2} U^{\varphi}[f_n] W[f_n]=\sup _{\varphi, \psi \in \bar{S}_{[0,1]}} \int_{[0,1]^2}U^{\varphi}[f_n] W^\psi[f_n]\\
& =\sup \left\{ \int_{[0,1]^2}U_0[f_n]W_0[f_n] :\left(\exists \varphi, \psi \in \bar{S}_{[0,1]}\right) U=U_0^{\varphi}, W=W_0^\psi\right\} .
\end{aligned}\end{equation}

Therefore, the $\cF-$overlay functional $\cC_{\cF}$ is invariant under weak isomorphism. Thus, we can consider the $\cF-$overlay functional as a functional on the space  $\UKernel\times \UKernel$ (or one of its subsets as $\UGraphon\times \UGraphon$). 

The following result gives us control over the $\cF-$overlay functional $\cC_{\cF}$ in terms of a finite approximation.

\begin{lemma}\label{boundOverlayFunctionsApproxFinite}
 Let $W,W^{\prime}$ measure valued graphons and $q=\sup\|W\|_{TV}\sup\|W^{\prime}\|_{TV}.$ Then we have 
 $$ \sup _{\varphi \in S_{[0,1]}}\sum^{N}_{n=1}2^{-n}\left\langle W[f_n], W^{\prime \varphi}[f_n]\right\rangle_{2,\R}-\frac{q}{N}\leq   \mathcal{C}_{\mathcal{F}}\left(W, W^{\prime}\right)   \leq  \sup _{\varphi \in S_{[0,1]}}\sum^{N}_{n=1}2^{-n}\left\langle W[f_n], W^{\prime \varphi}[f_n]\right\rangle_{2,\R}+\frac{q}{N}$$
\end{lemma}

\proof
By definition of convergence determining sequence we have that for $n \in \mathbb{N}, f_n$ takes values in $[0,1]$, and, therefore, $W[f_n]$ takes values in $[-\sup\|W\|_{TV},\sup\|W\|_{TV}]$ and $W^{\prime \varphi}[f_n]
$ (for any $\varphi$ measure-preserving and bijective) takes values in $[-\sup\|W^{\prime}\|_{TV},\sup\|W^{\prime}\|_{TV}]$ .%Remember also that $f_0=\mathbbm{1}$, and thus $U\left[f_0\right]-W\left[f_0\right] \equiv 0$. 
Therefore, $$\left\langle W[f_n], W^{\prime \varphi}[f_n]\right\rangle_{2,\R}=\int_{[0,1]^2} W[f_n] W^{\prime \varphi}[f_n]\in [-q,q].
$$
Upper bounding (respectively lower bounding) $\left\langle W[f_n], W^{\prime \varphi}[f_n]\right\rangle_{2,\R}$ by $q$ (respectively $-q$) for indices $n>N$ in the sum in the definition of $\mathcal{C}_{2,\mathcal{F}}\left(W, W^{\prime}\right) $ we get the result.

\endproof
\begin{remark}
We observe that if $W$ and $W^{\prime}$ are probability graphons the constant $q$ in Lemma \ref{boundOverlayFunctionsApproxFinite} is just $1$.   
\end{remark}

Moreover, by Lemma \ref{ProprLinearProbGraphons}, it follows

\begin{equation}
\begin{aligned}\label{identityFOverlay}
& \sup _{\varphi \in S_{[0,1]}}\sum^{N}_{n=1}2^{-n}\left\langle W[f_n], W^{\prime \varphi}[f_n]\right\rangle_{2,\R}= \sup _{\varphi \in S_{[0,1]}}\int_{[0,1]^2} W(\sum^{N}_{n=1} f_n 2^{-n}W^{\prime \varphi}[f_n])d x d y=\\
&\mathcal{C}(W,\sum^{N}_{n=1} f_n 2^{-n}W^{\prime}[f_n]).
 \end{aligned}
\end{equation}

From this identity, we obtain the following relationship between the overlay functional $\cC$ and the $\cF-$overlay functional $\cC_{\cF}.$
\begin{corollary}
     Let $W,W^{\prime}$ measure valued graphons and $q=\sup\|W\|_{TV}\sup\|W^{\prime}\|_{TV}.$ Then we have 
 $$ \mathcal{C}(W,\sum^{N}_{n=1} f_n 2^{-n}W^{\prime}[f_n])-\frac{q}{N}\leq   \mathcal{C}_{\mathcal{F}}\left(W, W^{\prime}\right)   \leq  \mathcal{C}(W,\sum^{N}_{n=1} f_n 2^{-n}W^{\prime}[f_n])+\frac{q}{N}$$
\end{corollary}
\proof
The proof follows directly by substituting the identity \eqref{identityFOverlay} in the statement of Lemma \ref{boundOverlayFunctionsApproxFinite}.
\endproof

To conclude this section we prove an explicit identity linking the overlay functional $\cC$ and the $\cF-$overlay functional $\cC_{\cF}.$ 
\begin{lemma}\label{Eq:FOverlayCOverlayInfty}
    Let $W,W^{\prime}$ probability graphons and $(f_n)_n=\cF$ a convergence determining sequence. Then we have 
 $$ \mathcal{C}(W,\sum^{\infty}_{n=1} f_n 2^{-n}W^{\prime}[f_n])=  \mathcal{C}_{\mathcal{F}}\left(W, W^{\prime}\right)   $$
\end{lemma}
\proof
Let's define with $\widehat{W}(\mathrm{d} x ,\mathrm{d} y,\mathrm{d} z)$ the probability measure on $[0,1]^2\times \Space,$
$$
\widehat{W}(\mathrm{d} x ,\mathrm{d} y,\mathrm{d} z)=W(x,y,\mathrm{d}z)\mathrm{d} x \mathrm{d} y
$$
By monotone convergence (the function $f_nW^{\prime }[f_n]$ is non-negative as $W^{\prime}$ is a probability graphon and $f_n$ is $[0,1]-$valued as an element of a convergence determining sequence) we have
$$
\begin{aligned}
&\sum^{\infty}_{n=1}2^{-n}\left\langle W[f_n], W^{\prime \varphi}[f_n]\right\rangle_{2,\R}=\sum^{\infty}_{n=1}  2^{-n}\int_{[0,1]^2} \int_{\Space}f_n(z)W^{\prime }[f_n](x,y)W(x,y,\mathrm{d}z)\mathrm{d} x \mathrm{d} y\\
&\sum^{\infty}_{n=1}  2^{-n}\int_{[0,1]^2} \int_{\Space}f_n(z)W^{\prime }[f_n](x,y)\widehat{W}(\mathrm{d} x ,\mathrm{d} y,\mathrm{d} z)=\int_{[0,1]^2} \int_{\Space}\sum^{\infty}_{n=1}  2^{-n}f_n(z)W^{\prime }[f_n](x,y)\widehat{W}(\mathrm{d} x ,\mathrm{d} y,\mathrm{d} z)\\
&=\int_{[0,1]^2} \int_{\Space}\sum^{\infty}_{n=1}  2^{-n}f_n(z)W^{\prime }[f_n](x,y)W(x,y,\mathrm{d}z)\mathrm{d} x \mathrm{d} y=\int_{[0,1]^2} W(\sum^{\infty}_{n=1} f_n 2^{-n}W[f_n])\mathrm{d} x \mathrm{d} y.
\end{aligned}
$$
We could also have used dominated convergence alternatively. This concludes the proof because the supremum does preserve equality.
\endproof

\subsection{Quotient sets of probability graphons}

In this section, we will study another invariant for probability graphons: quotient sets.

Let $W\in\UKernel$ be a finite signed measure-valued kernel and let $\mathcal{P}=\left\{S_1, \ldots, S_k\right\}$ be a partition of $[0,1]$ in $k$ measurable sets for $k\geq 1$, we define the quotient graph (or simply quotient) $W / \mathcal{P}$  as the measure edge-decorated and vertex-weighted graph on $[k]$, with node weights $\alpha_i(W / \mathcal{P})=\lambda\left(S_i\right)$ (where $\lambda$ denotes the Lebesgue measure on $[0,1])$ and as measure decoration of the edge $e=\{i,j\}$ the measure
$$
\beta_{i j}(W / \mathcal{P})=\frac{1}{\lambda\left(S_i\right) \lambda\left(S_j\right)} \int_{S_i \times S_j} W .
$$
To study right-convergence for probability graphons, we will consider all quotient graphs of a given signed measure-valued kernel at the same time. For a measure-valued kernel $W$ and probability distribution $\mathbf{a}$ on $[K]$, we denote by $\mathcal{Q}_{\mathbf{a}}(W)$ the set of quotients $L=W / \mathcal{P}$ with $\alpha(L)=\mathbf{a}.$ We denote by $\mathcal{Q}_k(W)$ the set of quotients for all partitions in $k$ sets.

We consider quotient sets as subsets of the space $$\R^k\times\cP(\Space)^{{k} \choose {2}}.$$
In particular, every quotient graph $W / \mathcal{P}$ in the quotient set $\mathcal{Q}_k(W)$ will be considered as an element of $\R^k\times\cP(\Space)^{{k} \choose {2}},$ where $\alpha=(\alpha_1(W / \mathcal{P}),\ldots,\alpha_k(W / \mathcal{P}))\in \R^k$ and $\beta(W / \mathcal{P})=(\beta_{i j}(W / \mathcal{P}))_{i,j\in [k]}\in \cP(\Space)^{{k} \choose {2}}.$

We observe that quotient sets can express the overlay functional when the second of the kernels involved is a $\CbFunct-$graph, as the following identity shows: 
\begin{equation}\label{OverlayQuotient}
    \mathcal{C}(W, H^{\beta})=\max _{L \in \mathcal{Q}_{\mathbf{a}}(W)} \sum_{i, j \in[q]} \alpha_i(L) \alpha_j(L) \beta_{i j}(L) (\beta_{i j}(H) )\text {. }
\end{equation}
We now define distances between graphs with weighted vertices and measure decorated edges that we will use to compare quotient graphs. 
For two vertex weighted and edge measure decorated graphs $G_{\tilde{\alpha}}^{\tilde{\beta}}$ and $H_{\alpha}^{\beta}$ (where $\tilde \alpha $ and $\alpha$ denote the two vertex weights vectors and $\beta$ and $\tilde\beta$ represent the measure decorations) on the vertex set $[k]$ we define the $d_1$ and $d_{\square}$ distance as 

$$
d_1\left(G^{\tilde\beta}, H^{\beta}\right)=\sum_{i\in[k]}|\alpha_i(G)-\alpha_i(H)|+\sum_{i, j \in[k]}d_{\mathcal{LP}}(\alpha_i(G)\alpha_j(G)\tilde{\beta}_{ij}(G),\alpha_i(H)\alpha_j(H)\beta_{ij}(H),\cdot))
%\left|\int_{S_i \times S_j} W-\int_{S_i^{\prime} \times S_j^{\prime}} W\right| .
$$
and 
\begin{equation}\label{Eqdsquare}
    d_{\square}\left(G^{\tilde \beta},H^{\beta}\right)=\sum_{i\in[k]}|\alpha_i(G)-\alpha_i(H)|+\sup_{S,T\subset [k]}d_{\mathcal{LP}}(\sum_{i\in S, j \in T}\alpha_i(G)\alpha_j(G)\tilde{\beta}_{ij}(G),\sum_{i\in S, j \in T}\alpha_i(H)\alpha_j(H)\beta_{ij}(H)).
\end{equation}

We will consider $d_1$ and $d_{\square}$ as distances on the space $\R^k\times\cP(\Space)^{{k} \choose {2}}.$ It is easy to see that $d_1$ and $d_{\square}$ are distances as product distances of distances.

By the definition of quotient graphs $L\in \mathcal{Q}_{\mathbf{a}}(W)\subset \mathcal{Q}_{k}(W) $ and $L^{\prime}\in \mathcal{Q}_{\mathbf{a}^{\prime}}(W) \subset\mathcal{Q}_k(U)$ we have
$$
d_1\left(L, L^{\prime}\right)=\left\|\mathbf{a}-\mathbf{a}^{\prime}\right\|_1+\sum_{i, j \in[k]}d_{\mathcal{LP}}(W(S_i \times S_j,\cdot),U(S_i^{\prime} \times S_j^{\prime},\cdot))
%\left|\int_{S_i \times S_j} W-\int_{S_i^{\prime} \times S_j^{\prime}} W\right| .
$$
and 
$$
d_{\square}\left(L, L^{\prime}\right)=\left\|\mathbf{a}-\mathbf{a}^{\prime}\right\|_1+\sup_{S,T\subset [k]}d_{\mathcal{LP}}(\sum_{i\in S, j \in T}W(S_i \times S_j,\cdot),\sum_{i\in S, j \in T}U(S_i^{\prime} \times S_j^{\prime},\cdot)).
$$

for suitable partitions $S_1,\ldots,S_k\subset [0,1]$ and $S^{\prime}_1,\ldots,S^{\prime}_k\subset [0,1]$ and where  $\left\|\mathbf{b}\right\|_1=\sum_{i\in[k]}|b_i|$ for any vector $\mathbf{b}\in \R^k$. 

We now show explicit bounds between $d_1$ and $d_{\square}.$

\begin{lemma}\label{Lemma_d1_dsquare} For $L\in \mathcal{Q}_k(W) $ and $L^{\prime}\in \mathcal{Q}_k(U)$ we have
$$
\frac{1}{k^2}d_1\left(L, L^{\prime}\right)\leq d_{\square}\left(L, L^{\prime}\right)\leq k^2d_1\left(L, L^{\prime}\right)
$$
\end{lemma}
\proof
We first prove $d_{\square}\left(L, L^{\prime}\right)\leq k^2d_1\left(L, L^{\prime}\right).$ This inequality follows from

$$\begin{aligned}
    d_{\square}\left(L, L^{\prime}\right)=&\left\|\mathbf{a}-\mathbf{a}^{\prime}\right\|_1+\sup_{S,T}d_{\mathcal{LP}}(\sum_{i\in S, j \in T}W(S_i \times S_j,\cdot),\sum_{i\in S, j \in T}U(S_i^{\prime} \times S_j^{\prime},\cdot))\leq\\
&   \left\|\mathbf{a}-\mathbf{a}^{\prime}\right\|_1+\sup_{S,T}|S||T|d_{\mathcal{LP}}\left(\frac{\sum_{i\in S, j \in T}W(S_i \times S_j,\cdot)}{|S||T|},\frac{\sum_{i\in S, j \in T}U(S_i^{\prime} \times S_j^{\prime},\cdot)}{|S||T|}\right)\leq \\
&\left\|\mathbf{a}-\mathbf{a}^{\prime}\right\|_1+\sup_{S,T}|S||T|\max_{i\in S, j \in T}d_{\mathcal{LP}}(W(S_i \times S_j,\cdot),U(S_i^{\prime} \times S_j^{\prime},\cdot))\leq  \\
& \left\|\mathbf{a}-\mathbf{a}^{\prime}\right\|_1+k^2\max_{i,j \in [k]}d_{\mathcal{LP}}(W(S_i \times S_j,\cdot),U(S_i^{\prime} \times S_j^{\prime},\cdot))\leq \\ 
    &    k^2d_1\left(L, L^{\prime}\right),
\end{aligned}$$

where in the first inequality we used Lemma \ref{LemmaIneqScalingProkhorov} and in the second inequality we used Lemma \ref{LemmaQuasi-convProhorov}.

Now the inequality $d_{1}\left(L, L^{\prime}\right)\leq k^2d_{\square}\left(L, L^{\prime}\right)$ follows from the estimate

$$
\begin{aligned}
d_1\left(L, L^{\prime}\right)=&
\left\|\mathbf{a}-\mathbf{a}^{\prime}\right\|_1+\sum_{i, j \in[q]}d_{\mathcal{LP}}(W(S_i \times S_j,\cdot),U(S_i^{\prime} \times S_j^{\prime},\cdot)) \leq \\
& \left\|\mathbf{a}-\mathbf{a}^{\prime}\right\|_1+k^2 \max_{i,j\in[k]}d_{\mathcal{LP}}(W(S_i \times S_j,\cdot),U(S_i^{\prime} \times S_j^{\prime},\cdot)) \leq \\
& k^2d_{\square}\left(L, L^{\prime}\right).
\end{aligned}
$$
\endproof

\begin{remark}
We proved Lemma \ref{Lemma_d1_dsquare} for quotients but it is easy to see that it holds for any measure edge decorated graph with weighted vertices.    
\end{remark}
To study probability graphons convergence we want to compare quotient sets, which are sets of quotient graphs. The Hausdorff metric is a natural choice to compare sets in a metric space.

\begin{definition}\label{DefHausdorffDist}
Let $(M, d)$ be a metric space. For each pair of non-empty subsets $X \subset M$ and $Y \subset M$, the \emph{Hausdorff distance} between $X$ and $Y$ is defined as
$$
d^{\text {Haus }}(X, Y):=\max \left\{\sup _{x \in X} d(x, Y), \sup _{y \in Y} d(X, y)\right\},
$$

 where $d(a, B)=\inf _{b \in B} d(a, b)$ quantifies the distance from a point $a \in M$ to the subset $B \subseteq M$.
\end{definition}

As anticipated, we will consider quotient graphs as elements of the space $M=\R^k\times\cP(\Space)^{{k} \choose {2}}$ equipped with the metric $d_{\square}$ or $d_1.$ To differentiate the Hausdorff metric for these two different choices of distance on $M$ we will denote by $d_{\square}^{\text {Haus }}$ the Hausdorff distance when we equip $M$ with the metric $d_{\square}$ and $d^{\text {Haus }}_1$ when we equip $M$ with the metric $d_{1}.$  

The next lemma gives a bound on the Hausdorff distance $d_1^{\text {Haus }}$ between quotient sets  $\mathcal{Q}_{\mathbf{a}}(W), \mathcal{Q}_{\mathbf{a}^{\prime}}(W)$ for the same probability graphon $W$ in terms of the probability distributions $\mathbf{a}$ and $\mathbf{a}^{\prime}.$

\begin{lemma}\label{LemmaBoundQuotients3L1}
For any $W \in \mathcal{W}_1$ and any two probability distributions  $\mathbf{a},\mathbf{a }^{\prime}$ on $[k]$, we have $$d_1^{\text {Haus }}\left(\mathcal{Q}_{\mathbf{a}}(W), \mathcal{Q}_{\mathbf{a}^{\prime}}(W)\right) \leq 3\left\|\mathbf{a}-\mathbf{a}^{\prime}\right\|_1.$$
\end{lemma} 

\proof Let $L=W / \mathcal{P} \in \mathcal{Q}_{\mathbf{a}}(W)$, where $\mathcal{P}=\left\{S_1, \ldots, S_k\right\} \in \Pi(\mathbf{a})$. We can easily construct a partition $\mathcal{P}^{\prime}=\left\{S_1^{\prime}, \ldots, S_k^{\prime}\right\} \in \Pi\left(\mathbf{a}^{\prime}\right)$ such that either $S_i \subseteq S_i^{\prime}$ or $S_i^{\prime} \subseteq S_i$ for every $i$. Let $L^{\prime}=W / \mathcal{P}^{\prime} \in \mathcal{Q}_{\mathbf{a}^{\prime}}(W)$. By definition,
$$
d_1\left(L, L^{\prime}\right)=\left\|\mathbf{a}-\mathbf{a}^{\prime}\right\|_1+\sum_{i, j \in[k]}d_{\mathcal{LP}}(W(S_i \times S_j,\cdot),W(S_i^{\prime} \times S_j^{\prime},\cdot))
%\left|\int_{S_i \times S_j} W-\int_{S_i^{\prime} \times S_j^{\prime}} W\right| .
$$
Moreover, by Lemma \ref{LemmBoundProhTotVar} we obtain
%\left|\int_{S_i \times S_j} W-\int_{S_i^{\prime} \times S_j^{\prime}} W\right| & \leq \lambda\left(\left(S_i \times S_j\right) \triangle\left(S_i^{\prime} \times S_j^{\prime}\right)\right) \\
$$
\begin{aligned}
&d_{\mathcal{LP}}(W(S_i \times S_j,\cdot),W(S_i^{\prime} \times S_j^{\prime},\cdot))\leq 
\\
&\|W(S_i \times S_j,\cdot)-W(S_i^{\prime} \times S_j^{\prime},\cdot)\|_{TV}\leq \\ 
&\leq \left\|\int_{S_i \times S_j\setminus S_i^{\prime} \times S_j^{\prime}} W-\int_{S_i^{\prime} \times S_j^{\prime}\setminus S_i \times S_j} W\right \|_{TV} \\
& \leq \sup \|W\|_{TV} \lambda\left(\left(S_i \times S_j \right) \triangle\left(S_i^{\prime} \times S_j^{\prime}\right)\right) \\
& \leq \sup \|W\|_{TV} (\left|a_i-a_i^{\prime}\right| \max \left(a_j, a_j^{\prime}\right)+\left|a_j-a_j^{\prime}\right| \min \left(a_i, a_i^{\prime}\right)) ,
\end{aligned}
$$
where $\triangle$ represents the symmetric set difference.
%{\color{blue} for previous inequality see also end of page 15 of \cite{abraham2023probabilitygraphons}}
Summing over all $i$ and $j$, we get
$$
\begin{aligned}
d_1\left(L, L^{\prime}\right) & \leq\left\|\mathbf{a}-\mathbf{a}^{\prime}\right\|_1\left(1+\sup\|W\|_{TV}(\sum_j \max \left(a_j, a_j^{\prime}\right)+\sum_i \min \left(a_i, a_i^{\prime}\right))\right) \\
& =\left\|\mathbf{a}-\mathbf{a}^{\prime}\right\|_1\left(1+\sup\|W\|_{TV}(\sum_j\left(a_j+a_j^{\prime})\right)\right)=(1+2\sup\|W\|_{TV})\left\|\mathbf{a}-\mathbf{a}^{\prime}\right\|_1 .
\end{aligned}
$$

This proves the lemma.
\endproof

We will need a couple of lemmas about the Hausdorff distance of quotient sets of different finite signed measure-valued kernels. We start showing that for a sequence of measure-valued kernels $(U_n)_n$ the convergence in unlabelled cut-metric in $\delta_{\square}$ implies convergence of the quotients sets ${Q}_k(U_n)$ in Hausdorff metric $d_{\square}^{\text {Haus }}.$ This follows from the following lemma.

\begin{lemma}\label{LemmMetricconvImplQuotient}
For any two finite signed measure-valued kernels  $U,W\in \Kernel$ and any integer $k \geq 1$, the following inequality holds:
$$
d_{\square}^{\text {Haus }}\left(\mathcal{Q}_k(U), \mathcal{Q}_k(W)\right) \leq \delta_{\square}(U, W) .
$$
\end{lemma} \proof Observe the invariance of the quotient sets under weak isomorphisms. Thus we consider, without loss of generality, the case in which $U$ and $W$ are optimally overlayed, i.e.\ that $\delta_{\square}(U, W)=d_{\square}\left(U, W\right)$. From Lemma 4.11 in \cite{abraham2023probabilitygraphons} %(Exercise 9.17 ) 
we have that $d_{\square}(U / \mathcal{P}, W / \mathcal{P}) \leq d_{\square}\left(U, W\right)$ for any measurable partition $\mathcal{P}$ of $[0,1]$ in $k$ sets. By the definition of Hausdorff distance, we obtain
$$
d_{\square}^{\text {Haus }}\left(\mathcal{Q}_k(U), \mathcal{Q}_k(W)\right) \leq  d_{\square}\left(U, W\right).
$$
This concludes the proof because quotient sets are invariant up to weak isomorphism.
\endproof

Another result that we will need to compare quotient sets later is the following lemma. 

\begin{lemma}\label{LemmaineqQuot}
For any two finite signed measure-valued kernels  $U,W\in \Kernel$ and any integer $k \geq 1$, the following inequalities hold:
$$
d_{\square}^{\text {Haus }}\left(\mathcal{Q}_k(U), \mathcal{Q}_k(W)\right) \leq \sup _{\mathbf{a}} d_{\square}^{\text {Haus }}\left(\mathcal{Q}_{\mathbf{a}}(U), \mathcal{Q}_{\mathbf{a}}(W)\right) \leq 4 d_{\square}^{\text {Haus }}\left(\mathcal{Q}_k(U), \mathcal{Q}_k(W)\right),
$$
where the supremum in the central term is taken over all $\mathbf{a}$ probability distributions on $[k].$
\end{lemma} \proof Let's start showing the first inequality. Let $H \in \mathcal{Q}_k(U)$, then $H \in \mathcal{Q}_{\mathbf{b}}(U)$ for the distribution $\mathbf{b}=\alpha(H).$ Hence
$$
\begin{aligned}
d_{\square}\left(H, \mathcal{Q}_k(W)\right) & \leq d_{\square}\left(H, \mathcal{Q}_{\mathbf{b}}(W)\right) \leq d_{\square}^{\text {Haus }}\left(\mathcal{Q}_{\mathbf{b}}(U), \mathcal{Q}_{\mathbf{b}}(W)\right) \\
& \leq \sup _{\mathbf{a}} d_{\square}^{\text {Haus }}\left(\mathcal{Q}_{\mathbf{a}}(U), \mathcal{Q}_{\mathbf{a}}(W)\right) .
\end{aligned}
$$
Since the previous inequality holds for every $H \in \mathcal{Q}_k(U)$ and, in the same way, for every element of $\mathcal{Q}_k(W)$, the result follows from the definition of the Hausdorff distance.

 Let's show now the second inequality. Let's $H \in \mathcal{Q}_{\mathbf{a}}(U)$ and $\mathbf{a}$ be any probability distribution on $[k].$  By the definition of Hausdorff metric, for every $\varepsilon>0$, there exists a quotient graph $L \in \mathcal{Q}_k(W)$ such that $d_{\square}(H, L) \leq d_{\square}^{\text {Haus }}\left(\mathcal{Q}_k(U), \mathcal{Q}_k(W)\right)+\varepsilon$. By the definition of $d_{\square}(H, L),$ for $L \in \mathcal{Q}_{\mathbf{b}}(W),$ we have $\|\mathbf{a}-\mathbf{b}\|_1 \leq d_{\square}(H, L)$. Moreover, by Lemma \ref{LemmaBoundQuotients3L1}, there exists a quotient graph $L^{\prime} \in \mathcal{Q}_{\mathbf{a}}(W)$ such that $d_{\square}\left(L, L^{\prime}\right) \leq k^2d_1\left(L, L^{\prime}\right) \leq k^2(3\|\mathbf{a}-\mathbf{b}\|_1+\varepsilon )\leq 3k^2 d_{\square}(H, L)+k^2\varepsilon$. Hence we obtain
$$
d_{\square}\left(H, L^{\prime}\right) \leq k^2d_1(H, L)+d_{\square}\left(L, L^{\prime}\right) \leq 4 k^2d_{\square}(H, L)+k^2\varepsilon \leq 4 k^2d_{\square}^{\text {Haus }}\left(\mathcal{Q}_k(U), \mathcal{Q}_k(W)\right)+(k^2+1) \varepsilon.
$$
This shows the lemma as $\varepsilon>0$ was chosen arbitrarily.
\endproof

\subsection{A random graph model}
This short section presents a random $\cC-$graph model, where $\cC$ is a finite set of continuous and bounded functions $\cC=\{f_1,\ldots,f_N\}\subset\CbFunct,$ obtained by a probability graphon $U=\sum^N_{i=1}f_iw_i,$ where $w_1,\ldots,w_N$ are real-valued graphons.

Let $f_1,\ldots,f_N\in \CbFunct,$ continuous and bounded (taking values in [0,1]) functions, $w_1,\ldots,w_N$ real-valued graphons (taking values in $[0,1]$) such that $\sum^N_{i=1}w_i(x,y)=1$ for almost every $(x,y)\in [0,1]^2.$ 

\begin{remark}
Let's assume that $\sum^N_{i=1}w_i(x,y)\leq 1$ for each $x,y\in [0,1].$ In this case, one can enrich $\cC$ with the function $f_{N+1}=0$ and define the real-valued graphon $w_{N+1}=\mathbbm{1}_{[0,1]^2}-\sum^N_{i=1}w_i.$  For this reason, the requirement $\sum^N_{i=1}w_i(x,y)=1$ can always be relaxed to $\sum^N_{i=1}w_i(x,y)\leq 1.$
\end{remark}
Let $U$ be the $\CbFunct-$valued kernel $U=\sum^N_{i=1}f_iw_i$, and $x=\left(x_1, \ldots, x_n\right)\in [0,1]^n,$  for $n \in \mathbb{N}.$ We define the probability measures (directed) edge decorated graph $\mathbb{H}(x, W)$ as the complete graph with vertex set $[n]=\{1, \ldots, n\}$, and with each (directed) edge $(j, k)$ decorated by the probability measure $\P_{jk}$ on $\{f_1,\ldots,f_N\}$ defined as $ \P_{jk}(\{f_i\})=w_i(x_j,x_k)$ (this already characterize the measure completely because it is a finite measure).

We can define from $\mathbb{H}(x, U)$ a random (directed) $\{f_1,\ldots,f_N\}-$graph $\mathbb{G}(x,U)$ as the complete graph whose vertex set is $[n]$, and with each edge $(j, k)$ having a random decoration $\beta_{jk}\in \{f_1\ldots,f_N\}\subset\CbFunct$ distributed according to the probability distribution $\P_{jk}$, all the random decorations being independent of each other. 

In particular, we will be interested in the special case when the sequence $x\in [0,1]^n$ is chosen uniformly at random: Let $X=\left(X_i\right)_{1 \leq i \leq n}$ be random variables where $X_i$ are independent and uniformly distributed on $[0,1].$ For this special case, we will denote $\mathbb{G}(n, U)=\mathbb{G}(X, U)$, which conditioning on $X=x$ is equal to $\mathbb{G}(x, U)$. We will call the random $\{f_1,\ldots,f_N\}-$graph $\mathbb{G}(n, U)$ the $U$-\emph{random graph} on $[n]$.

Hence, for a realization $G_{\omega}$ of $\mathbb{G}(n,U),$ which is a $\{f_1,\ldots,f_N\}-$graph, we can write the formal sum $G_{\omega}=\sum^N_{i=1}f_iH_{i,\omega}^{(n)},$ where $H_i^{(n)}$ are simple graphs such that $e\in E(H_{i,\omega}^{(n)}) $  if and only if $e\notin E(H_{j,\omega}^{(n)}) $ for every $j\in [n]$ such that $j\neq i.$

\begin{remark}
The $U$-random graph $\mathbb{G}(n, U)$ is in general directed. There can be $\beta_{jk}\neq \beta_{kj}$ also when $U$ is symmetric, i.e.\ $U(x,y)=U(y,x)$ for a.e.\ $(x,y)\in [0,1]^2.$ However, one can consider a symmetrized version of $\mathbb{G}(n, U)$ that we call $\mathbb{G}^{\text{sym}}(n, U)$. Let $\beta_{jk}$ for $j,k\in [n]$ be the random decoration of $\mathbb{G}(n, U).$ The random graph $\mathbb{G}^{\text{sym}}(n, U)$ is the complete graph on the vertex set $[n]$ with random edge decoration $\Tilde{\beta}_{jk}$  where $\Tilde{\beta}_{jk}=\beta_{jk}$ if $j\geq k$ and $\Tilde{\beta}_{jk}=\beta_{kj}$ if $j\leq k.$
\end{remark}

We show now that these random graph models can approximate the probability graphons they are sampled from.

For a simple graph $H$ on the vertex set $V=\{v_1,\ldots,v_n\}$ we denote by $w_H$ the real-valued graphon $$
w_H(x,y)=\begin{cases}
    1 & if \ (x,y)\in I_j=[j-1/n,j/n] \\
    0 & elsewhere.
\end{cases}
$$
For a $\CbFunct-$graph $G=\sum^N_{j=1}f_jH_j$ we define $$w_G=\sum^N_{j=1}f_jw_{H_j}.$$
\begin{lemma}[Convergence of sampled $\CbFunct-$graphs]\label{LemmSampledCbGraphsConvergeToKernel}
Let $f_1,\ldots,f_N\in \CbFunct,$ continuous and bounded (taking values in [0,1]) functions, $w_1,\ldots,w_N$ real-valued graphons (taking values in $[0,1]$) such that $\sum^N_{i=1}w_i(x,y)=1$ for almost every $(x,y)\in [0,1]^2.$ Let $U=\sum^N_{i=1}f_iw_i$ and consider $\mathbb{G}(n,U)=\sum^N_{i=1}f_iH_i^{(n)},$ the random $\cF-$graph on  the vertex set $[n]$ sampled from $U.$

For the sequence of $\cF-$valued kernels $U_n=\sum^N_{i=1}f_iw_{H^{(n)}_i}$ (obtained by the graphs $\mathbb{G}(n,U)$) we have (recall \eqref{Eq:CutmetricFiniteCb})
\begin{equation*}
\delta_{\square,N}(U,U_n)=\inf_{\varphi}\sum^N_{i=1}\|w_i-w^{\varphi}_{H^{(n)}_i}\|_{\square,\R}\rightarrow 0
\end{equation*} almost surely.
\end{lemma}
\proof
One can identify $f_1,\ldots,f_N\in \CbFunct$ with the Dirac measures $\delta_1,\ldots,\delta_N$ (where we recall that $\delta_i$ denotes the Dirac measure centered at $i\in [N]$). Therefore, we can identify $U=\sum^N_{i=1}f_iw_i$ with the probability graphon $U_p=\sum^N_{i=1}\delta_iw_i.$ The lemma then follows from the result for probability graphons (Theorem 6.13 in \cite{abraham2023probabilitygraphons}) as we are considering the discrete topology/metric both on the set $\{f_1,\ldots,f_N\}$ and on the set  $\{\delta_1,\ldots,\delta_N\}$. 
%An $\{f_1,\ldots f_n\}-$graph can approximate $\CbFunct-$valued kernels of the form $U=\frac{1}{N}\sum^N_{n=1}f_nw_n$ coordinatewise for every $w_n$ in $\NcutR{}.$ %follows partitioning $[0,1]$ into $N$ equal intervals $I_i$, and therefore obtaining a partition of the unit square $[0,1]^2$  of squares $I_i\times I_j$ for $i,j\in [N].$ 
%We sample $x_1,\ldots ,x_k$ i.i.d.\ uniformly in $[0,1].$ For every element $(i,j)\in[k]\times [k]$ we sample an element of $n\in[N]$ uniformly at random. Then we add edge $(i,j)$ with decoration $f_n$ with probability $w_n(x_i,x_j).$
\endproof

\subsection{Probability graphon convergence from the right}

In this section, we finally prove the main result of this work Theorem \ref{ThmEquivalenceConvQuotientOverlay} and Theorem \ref{ThmEquivalenceConvQuotientOverlayVersion2}. 

We start proving a characterization of weak isomorphism for probability graphons in terms of the overlay functional $\cC.$

\begin{lemma}\label{EquivWeakIsoOverl}
Let $\cF$ be a convergence-determining sequence. Let $W,W^{\prime}$ be two probability graphons. The following are equivalent:

\begin{enumerate}
\item $W$ and $W^{\prime}$ are weakly isomorphic $(W\sim W^{\prime})$.
\item $\cC(W,U)=\cC(W^{\prime},U)$
for every (Bochner integrable) $\CbFunct-$valued kernel $U;$
\item $\cC(W,U)=\cC(W^{\prime},U)$
for every $\cF-$valued kernel $U;$
\item $\cC(W,H^{\beta})=\cC(W^{\prime},H^{\beta})$
for every $\cF-$graph $H^{\beta};$
\item  $\cC_{\cF}(W,Q)=\cC_{\cF}(W^{\prime},Q)\left(=\cC_{\cF}(Q,W)=\cC_{\cF}(Q,W^{\prime})\right)$ for every probability graphon $Q.$

\end{enumerate}    
\end{lemma}
\proof
(1) $\Rightarrow$ (2) follows from Corollary \ref{CorSequencesEqOv}. \newline
(2) $\Rightarrow$ (3) is trivial. \newline
(3) $\Rightarrow$ (4) is trivial by identity \eqref{eq:IdentOverlayGraphGraphon}. \newline
(4) $\Rightarrow$ (1) We prove the contrapositive. Let's assume that $W$ and $W^{\prime}$ are not weakly isomorphic. This means that $\delta_{\square,\mathcal{F}}(W,W^{\prime})>0.$ From Lemma 4.15 in \cite{abraham2023probabilitygraphons} we obtain that 

$$
0<\delta_{\square,\mathcal{F}}(W,W^{\prime})\leq \sqrt{2}\delta_{2,\mathcal{F}}(W,W^{\prime}).
$$
Therefore, there exists a probability graphon $U$ such that $\mathcal{C}(W, U) \neq \mathcal{C}\left(W^{\prime}, U\right).$ In fact, \eqref{EqOverlayFunctConvDetFam} implies
$$
\left(\mathcal{C}_{2,\mathcal{F}}\left(W^{\prime}, W^{\prime}\right)-\mathcal{C}_{2,\mathcal{F}}\left(W^{\prime}, W\right)\right)+\left(\mathcal{C}_{2,\mathcal{F}}(W, W)-\mathcal{C}_{2,\mathcal{F}}\left(W^{\prime}, W\right)\right)=\delta_{2,\mathcal{F}}\left(W^{\prime}, W\right)^2>0,
$$

and so either $\mathcal{C}_{2,\mathcal{F}}\left(W^{\prime}, W^{\prime}\right) \neq \mathcal{C}_{2,\mathcal{F}}\left(W, W^{\prime}\right)$ or $\mathcal{C}_{2,\mathcal{F}}\left(W, W^{\prime}\right) \neq \mathcal{C}_{2,\mathcal{F}}(W, W)$. We will assume without loss of generality that $\mathcal{C}_{2,\mathcal{F}}\left(W^{\prime}, W^{\prime}\right) \neq \mathcal{C}_{2,\mathcal{F}}\left(W, W^{\prime}\right),$ the other case works in the exact same way.%i.e.\

Using Lemma \ref{boundOverlayFunctionsApproxFinite} we get that there exists $N\in\N$ large enough such that 
$$
 \sup _{\varphi \in S_{[0,1]}}\sum^{N}_{n=1}2^{-n}\left\langle W^{\prime}[f_n], W^{\prime \varphi}[f_n]\right\rangle_{2,\R}\neq\sup _{\varphi \in S_{[0,1]}}\sum^{N}_{n=1}2^{-n}\left\langle W[f_n], W^{\prime \varphi}[f_n]\right\rangle_{2,\R}.
$$
Moreover, using Lemma \ref{ProprLinearProbGraphons} (or equation \eqref{identityFOverlay}), we obtain

$$
 \sup _{\varphi \in S_{[0,1]}}\sum^{N}_{n=1}2^{-n}\left\langle W[f_n], W^{\prime \varphi}[f_n]\right\rangle_{2,\R}= \sup _{\varphi \in S_{[0,1]}}\int_{[0,1]^2} W(\sum^{N}_{n=1} f_n 2^{-n}W^{\prime \varphi}[f_n])d x d y=\mathcal{C}(W,\sum^{N}_{n=1} f_n 2^{-n}W^{\prime}[f_n]).
$$
and similarly

$$
 \sup _{\varphi \in S_{[0,1]}}\sum^{N}_{n=1}2^{-n}\left\langle W^{\prime}[f_n], W^{\prime \varphi}[f_n]\right\rangle_{2,\R}%= \sup _{\varphi \in S_{[0,1]}}\int_{[0,1]^2} W(\sum^{N}_{n=1} f_n 2^{-n}W^{\prime \varphi}[f_n])d x d y
 =\mathcal{C}(W^{\prime},\sum^{N}_{n=1} f_n 2^{-n}W^{\prime}[f_n]).
$$
Therefore, 
$$
\mathcal{C}(W,\sum^{N}_{n=1} f_n 2^{-n}W^{\prime}[f_n])\neq \mathcal{C}(W^{\prime},\sum^{N}_{n=1} f_n 2^{-n}W^{\prime}[f_n]). 
$$
Using the homogeneity of the overlay functional we thus have
$$
\mathcal{C}(W,\sum^{N}_{n=1} \frac{1}{N}f_n 2^{-n}W^{\prime}[f_n])\neq \mathcal{C}(W^{\prime},\sum^{N}_{n=1} \frac{1}{N}f_n 2^{-n}W^{\prime}[f_n]). 
$$
We can now set $w_n=2^{-n}W^{\prime}[f_n]$ and $U=\frac{1}{N}\sum^N_{n=1}f_nw_n.$ We observe that $w_1,\ldots,w_n$ are real-valued graphons taking values in $[0,1]$ because $W^{\prime}$ is a probability graphon and the functions $f_1,\ldots,f_N$ take values in $[0,1]$ as they are elements of a convergent determining sequence. We have just shown that $U=\frac{1}{N}\sum^N_{n=1}f_nw_n$ is a $span(\cF)-$valued kernel such that $$
\mathcal{C}(W,U)\neq \mathcal{C}(W^{\prime},U). 
$$
%In particular, we have just proven that we can choose $U=\frac{1}{N}\sum^N_{n=1}f_nw_n$ where $f_1,\ldots,f_n\in \cF$ and $w_1,\ldots,w_n$ are real-valued graphons (i.e.\ they take values in $[0,1]$). This follows because $W^{\prime}\in \UGraphon$ and for every $n\in\N$ the function $f_n$ takes values in $[0,1]$ as element of a convergent determining sequence and we can rescale $U$ with a positive scalar keeping the inequality for "homogeneity". 
Furthermore, we can choose $w_i$ of the form $w_i=w_{H_i}$, where $H_i$ is a simple graph and $w_i,w_j$ have disjoint supports for $i\neq j$. This follows from Lemma \ref{LemmSampledCbGraphsConvergeToKernel} since $\{f_1,\ldots f_n\}-$graphs can approximate $\{f_1,\ldots f_n\}-$valued kernels of the form $U=\sum^N_{n=1}f_n(\frac{1}{N}w_n)$ in $\delta_{\square,N}$ (see equation \eqref{Eq:CutmetricFiniteCb} for the definition) and from the continuity of the overlay functional with respect to $\delta_{\square,N}$ for $\CbFunct-$valued kernels $U$ of the form $U=\sum^N_{i=1}f_iw_i$ where $f_i\in\CbFunct$ and $w_i$ are real-valued kernels, Lemma \ref{continuityContBoundGraphRealValCutDist}.\newline
This concludes the proof.

%We justify here the fact that simple $\{f_1,\ldots f_n\}-$graphs can approximate $\CbFunct-$valued kernels of the form $U=\frac{1}{N}\sum^N_{n=1}f_nw_n$ coordinatewise for every $w_n$ in $\NcutR{}.$ %follows partitioning $[0,1]$ into $N$ equal intervals $I_i$, and therefore obtaining a partition of the unit square $[0,1]^2$  of squares $I_i\times I_j$ for $i,j\in [N].$ 
%We sample $x_1,\ldots ,x_k$ i.i.d.\ uniformly in $[0,1].$ For every element $(i,j)\in[k]\times [k]$ we sample an element of $n\in[N]$ uniformly at random. Then we add edge $(i,j)$ with decoration $f_n$ with probability $w_n(x_i,x_j).$
(1) $\Rightarrow$ (5) Follows from identity \eqref{EqWeakIsoFOverlaFunc}. \newline
(5) $\Rightarrow$ (1) We observe that we already proved this implication in the proof of the implication (4) $\Rightarrow$ (1).
\endproof

%{We will also use the following lemma which is a direct corollary of the Uniform Boundedness Principle (Theorem 2.2 in \cite{brezis2010functional}).

%\color{blue} Corollary 2.3 in \cite{brezis2010functional}
%\begin{lemma}\label{LemmCorUnifBoundPrin}
%Let $X$ and $Y$ be Banach spaces. If $\left(T_n\right)$ is a sequence of bounded linear operators from $X$ to $Y$ such that for every $x \in X, \lim _{n \rightarrow \infty} T_n(x)$ exists, then $T(x)=\lim _{n \rightarrow \infty} T_n(x)$ is a bounded linear operator from $X$ to $Y$.
%\end{lemma}
%}

 We now state the first version of our main result.

\begin{theorem}\label{ThmEquivalenceConvQuotientOverlay}
For any tight sequence $\left(W_n\right)$ in $\UGraphon$, the following are equivalent:
\begin{enumerate}

\item  The sequence $\left(W_n\right)$ is convergent in the cut distance $\delta_{\square};$
\item  The sequence $\left(W_n\right)$ is convergent in the cut distance $\delta_{\square,\cF}$ for any choice of convergence determining sequence $\cF;$
\item The overlay functional values $\mathcal{C}\left(W_n, U\right)$ are convergent for every $\CbFunct-$valued kernel $U$ Bochner measurable;
\item The overlay functional values $\mathcal{C}\left(W_n, H^{\beta}\right)$ are convergent for every decorated $\CbFunct-$graph $H^{\beta}$;
%\item The overlay functional values $\mathcal{C}\left(W_n, U\right)$ are convergent for every $\overline{span(\cF)}-$valued kernel $U$;
%\item The overlay functional values $\mathcal{C}\left(W_n, H\right)$ are convergent for every decorated $\overline{span(\cF)}-$graph $H$;
%\item The overlay functional values $\mathcal{C}\left(W_n, U\right)$ are convergent for every $span(\cF)-$valued kernel $U$;
%\item The overlay functional values $\mathcal{C}\left(W_n, H\right)$ are convergent for every decorated $span(\cF)-$graph $H$;
\item The overlay functional values $\mathcal{C}\left(W_n, U\right)$ are convergent for every $\cF-$valued kernel $U$;
\item The overlay functional values $\mathcal{C}\left(W_n, H^{\beta}\right)$ are convergent for every decorated $\cF-$graph $H^{\beta}$;
\item The quotient sets $\mathcal{Q}_k\left(W_n\right)$ form a Cauchy sequence in the $d_{\square}^{\text {Haus }}$ Hausdorff metric for every $k \geq 1$.
\item The $\cF-$overlay functional values $\mathcal{C_{\cF}}\left(W_n, Q\right)=\mathcal{C_{\cF}}\left(Q,W_n\right)$ are convergent for every probability graphon $Q.$
\end{enumerate}
\end{theorem}
\proof

(1) if and only if (2) follows from Corollary 5.6 in \cite{abraham2023probabilitygraphons} (see Theorem \ref{cor:equiv-topo}). \newline
(2) $\Rightarrow$ (3) by Corollary \ref{CorSequencesEqOv}. \newline
(3) $\Rightarrow$ (4) is trivial by identity \eqref{eq:IdentOverlayGraphGraphon}\newline
(3) $\Rightarrow$ (5) is trivial.\newline
(5) $\Rightarrow$ (6) is trivial by identity \eqref{eq:IdentOverlayGraphGraphon}.\newline
(6) $\Rightarrow$ (1) We prove the contrapositive. Let $\left(W_n\right)$ be a sequence of probability graphons which is not convergent in the unlabelled cut distance $\delta_{\square,\mathcal{LP}}$. By the relative compactness of tight subsets of the probability graphon space (Theorem \ref{ThmRelCompact} which is Proposition 5.2 in \cite{abraham2023probabilitygraphons}), $\left(W_n\right)$ has two subsequences $\left(W_{n_i}\right)$ and $\left(W_{m_i}\right)$ converging to not weakly isomorphic probability graphons $W$ and $W^{\prime}$. 

From Lemma \ref{EquivWeakIsoOverl}, there exists an $\cF-$graph $H^{\beta}$ such that
$$
\mathcal{C}\left(W, H^{\beta}\right)\neq \mathcal{C}\left(W^{\prime}, H^{\beta}\right).$$
Therefore, it follows that the values $\mathcal{C}\left(W_n, H^{\beta}\right)$ cannot form a convergent sequence, contradicting statement (6), since $\mathcal{C}\left(W_{n_i}, H^{\beta}\right) \rightarrow \mathcal{C}(W, H^{\beta})$ and $\mathcal{C}\left(W_{m_i}, H^{\beta}\right) \rightarrow \mathcal{C}\left(W^{\prime}, H^{\beta}\right)$ and $$
\mathcal{C}\left(W, H^{\beta}\right)\neq \mathcal{C}\left(W^{\prime}, H^{\beta}\right).
$$\newline

This already proves that (1),(2),(3),(4),(5),(6).\newline

(1) $\Rightarrow$ (7) by Lemma \ref{LemmMetricconvImplQuotient}. \newline
(7) $\Rightarrow$ (6) Fix any $span(\cF)-$graph $H$ on $[k]$, and let $\mathbf{a}$ be the uniform distribution on $[k]$. Let $n, m \geq 1$, then from identity \eqref{OverlayQuotient} we have
$$
\mathcal{C}\left(W_n, H^{\beta}\right)=\sup _{L \in \mathcal{Q}_{\mathbf{a}}\left(W_n\right)} \frac{1}{k^2} \sum_{\substack{i, j \\ i j \in E(H)}} \beta_{i j}(L)(\beta_{i j}(H))=\max _{L \in \overline{\mathcal{Q}}_{\mathbf{a}}\left(W_n\right)} \frac{1}{k^2} \sum_{\substack{i, j \\ i j \in E(H)}} \beta_{i j}(L)(\beta_{i j}(H)) .
$$
where $\overline{\mathcal{Q}}_{\mathbf{a}}\left(W_n\right)$ is the closure of ${\mathcal{Q}}_{\mathbf{a}}\left(W_n\right)$ with respect to $d_{\square}$ (or equivalently $d_1$).
Let $L_n \in \overline{\mathcal{Q}}_{\mathbf{a}}\left(W_n\right)$ attain the maximum. By the definition of Hausdorff distance, there is an $L^{\prime} \in \overline{\mathcal{Q}}_{\mathbf{a}}\left(W_m\right)$ such that $d_{\square}\left(L_n, L^{\prime}\right) \leq d_{\square}^{\text {Haus }}\left(\mathcal{Q}_{\mathbf{a}}\left(W_n\right), \mathcal{Q}_{\mathbf{a}}\left(W_m\right)\right)$. The definition of the overlay functional $\mathcal{C}$ implies that $\mathcal{C}\left(W_m, H\right) \geq\left(1 / k^2\right) \sum_{i j \in E(H)} \beta_{i j}\left(L^{\prime}\right)(\beta_{i j}(H))$. Hence, for any $\varepsilon>0$ it exists $\delta>0$ small enough such that for $d_{\square}^{\text {Haus }}\left(\mathcal{Q}_k\left(W_n\right), \mathcal{Q}_k\left(W_n\right)\right)<\delta$ we obtain that $\mathcal{C}\left(W_n, H\right)  -\mathcal{C}\left(W_m, H\right)\leq \varepsilon $. This follows from 
$$
\begin{aligned}
\mathcal{C}\left(W_n, H^{\beta}\right) & -\mathcal{C}\left(W_m, H^{\beta}\right) \leq \frac{1}{k^2} \sum_{i j \in E(H)} \beta_{i j}\left(L_n\right)(\beta_{i j}(H))-\frac{1}{k^2} \sum_{\substack{i, j \\
i j \in E(H)}} \beta_{i j}\left(L^{\prime}\right)(\beta_{i j}(H)) \\
& \leq \frac{1}{k^2} \sum_{i, j\in [k]}\left|(\beta_{i j}\left(L_n\right)-\beta_{i j}\left(L^{\prime}\right))(\beta_{i j}(H))\right|\leq \varepsilon
\end{aligned}$$
by Portmanteau theorem (see for example Theorem 2.1 in \cite{billingsley1968convergence}). In fact, we can bound the Levy-Prokhorov distance $d_{\mathcal{LP}}$ (that is metrizing weak convergence) between $\frac{1}{k^2}\beta_{ij}(L_n)$ and $\frac{1}{k^2}\beta_{ij}(L^{\prime})$ in the following way 
$$\begin{aligned}
&d_{\mathcal{LP}}(\frac{1}{k^2}\beta_{ij}(L_n),\frac{1}{k^2}\beta_{ij}(L^{\prime})))\\
&= d_1\left(L_n, L^{\prime}\right)\\ &\leq k^2 d_{\square}\left(L_n, L^{\prime}\right) \\
& \leq k^2 d_{\square}^{\text {Haus }}\left(\mathcal{Q}_{\mathbf{a}}\left(W_n\right), \mathcal{Q}_{\mathbf{a}}\left(W_m\right)\right)\\
&\leq 4k^2 d_{\square}^{\text {Haus }}\left(\mathcal{Q}_k\left(W_n\right), \mathcal{Q}_k\left(W_m\right)\right),\end{aligned}
$$
where the second last inequality follows from Lemma \ref{LemmaineqQuot}.  
By Lemma \ref{LemmaIneqScalingProkhorov} we obtain $$\begin{aligned}
&d_{\mathcal{LP}}(\beta_{ij}(G),\beta_{ij}(H))\\
&\leq 4k^4 d_{\square}^{\text {Haus }}\left(\mathcal{Q}_k\left(W_n\right), \mathcal{Q}_k\left(W_m\right)\right)\end{aligned}
$$
and the right-hand side tends to 0 as $n, m \rightarrow \infty$ by hypothesis. This implies that
$$
\limsup _n\left(\mathcal{C}\left(W_n, H\right)-\mathcal{C}\left(W_m, H\right)\right) \leq 0 .
$$
Observing that the same argument works exchanging $n$ and $m,$ we obtain that $\left(\mathcal{C}\left(W_n, H\right): n=1,2, \ldots\right)$ is a Cauchy sequence.\newline

Therefore, we obtain that (1), (2), (3), (4), (5), (6), (7) are equivalent.

(3) $\Rightarrow$ (8) By Lemma \ref{Eq:FOverlayCOverlayInfty}. \newline
(8) $\Rightarrow$ (1) We prove the contrapositive. Let $\left(W_n\right)$ be a sequence of probability graphons which is not convergent in the unlabelled cut distance $\delta_{\square,\mathcal{LP}}$. By the relative compactness of tight subsets of the probability graphon space (Theorem \ref{ThmRelCompact} which is Proposition 5.2 in \cite{abraham2023probabilitygraphons}), $\left(W_n\right)$ has two subsequences $\left(W_{n_i}\right)$ and $\left(W_{m_i}\right)$ converging to not weakly isomorphic probability graphons $W$ and $W^{\prime}$. 

From Lemma \ref{EquivWeakIsoOverl}, there exists a probability graphon $Q$ such that
$$
\mathcal{C}_{\cF}\left(W, Q\right)\neq \mathcal{C}_{\cF}\left(W^{\prime}, Q\right).$$
Therefore, it follows that the values $\mathcal{C}_{\cF}\left(W_n, Q\right)$ cannot form a convergent sequence, contradicting statement (8), since $\mathcal{C}_{\cF}\left(W_{n_i}, Q\right) \rightarrow \mathcal{C}_{\cF}(W, Q)$ and $\mathcal{C}_{\cF}\left(W_{m_i}, Q\right) \rightarrow \mathcal{C}_{\cF}\left(W^{\prime}, Q\right).$ \newline
\endproof

\begin{remark}
We observe that the proof of the implication (6) $\Rightarrow$ (1) in the proof of Theorem \ref{ThmEquivalenceConvQuotientOverlay} is not effective.
%Some of the arguments in the proof of Theorem \ref{ThmEquivalenceConvQuotientOverlay}, most notably the proof of (6) $\Rightarrow$ (1), are not effective. 
One could try to prove an explicit inequality for (6) $\Rightarrow$ (1) as done in the case of real-valued graphons in \cite{borgs2011convergentAnnals}, but this is out of the scope of this work.
\end{remark}

\begin{remark}\label{RemarkMainThmMorGen}
   From conditions (4) and (6) it would be equivalent to assume the convergence of the sequence $\mathcal{C}\left(W_n, H^{\beta}\right)$ for every decorated $span(\cF)-$graph $H^{\beta}$ or for every decorated $\overline{span(\cF)}-$graph $H^{\beta}$.
From Lemma \ref{Lemma_d1_dsquare}, in (7) one could also use the $d_1^{\text {Haus }}$ Hausdorff metric as an alternative to the  $d_{\square}^{\text {Haus }}$ Hausdorff metric.    
    
From Lemma \ref{LemmaineqQuot} it follows that we could alternatively require in (7) the convergence of $\mathcal{Q}_{\mathbf{a}}\left(W_n\right)$ for every $k \geq 1$ and probability distribution $\mathbf{a}$ on $[k]$. 
 
It would be enough to require this only for the uniform distribution. We explain here the idea of the proof of this statement. Let's assume that  $\mathcal{Q}_{\mathbf{u}}\left(W_n\right)$ is a Cauchy sequence in $d_{\square}^{\text {Haus }}$ for every $k \geq 1$ where $\mathbf{u}$ is the uniform measure on $[k].$ Let's now consider the sequence $\mathcal{Q}_{\mathbf{a}}\left(W_n\right)$ for a fixed $s \geq 1$ and a fixed probability distribution $\mathbf{a}=(a_1,\ldots,a_s)$ on $[s].$ We therefore can select a $k\in \N$ big enough such that there exists a probability distribution $\mathbf{b}=(b_1,\ldots,b_s)$ on $[s]$ such that for every $i\in [s]$ we have $b_i=p_i/k$ for some $p_i\in [k]$ and $\|\mathbf{a}-\mathbf{b}\|_1<\varepsilon.$ Thus, from Lemma \ref{Lemma_d1_dsquare} and Lemma \ref{LemmaBoundQuotients3L1} we have that $
d_{\square}^{\text {Haus }}(\mathcal{Q}_{\mathbf{a}}\left(W_n\right),\mathcal{Q}_{\mathbf{b}}\left(W_n\right))<3s^2\varepsilon
$ for every $n\in N.$ Additionally, it follows from the definition of  $d_{\square},$ \eqref{Eqdsquare}, that
$$
d_{\square}^{\text {Haus }}(\mathcal{Q}_{\mathbf{b}}\left(W_n\right),\mathcal{Q}_{\mathbf{b}}\left(W_m\right))\leq d_{\square}^{\text {Haus }}(\mathcal{Q}_{\mathbf{u}}\left(W_n\right),\mathcal{Q}_{\mathbf{u}}\left(W_m\right)).
$$
and therefore also $\mathcal{Q}_{\mathbf{b}}\left(W_n\right)$ is a Cauchy sequence in $d_{\square}^{\text {Haus }}$.
Moreover, from the triangular inequality we obtain $$
d_{\square}^{\text {Haus }}(\mathcal{Q}_{\mathbf{a}}\left(W_n\right),\mathcal{Q}_{\mathbf{a}}\left(W_m\right))\leq2 s^2\varepsilon+d_{\square}^{\text {Haus }}(\mathcal{Q}_{\mathbf{b}}\left(W_n\right),\mathcal{Q}_{\mathbf{b}}\left(W_m\right)).
$$
Thus, also $\mathcal{Q}_{\mathbf{a}}\left(W_n\right)$ is a Cauchy sequence in $d_{\square}^{\text {Haus }}$ as we can choose $\varepsilon>0$ arbitrarily small. 
\end{remark}

We give also a slightly different version of the previous theorem here.

\begin{theorem}\label{ThmEquivalenceConvQuotientOverlayVersion2}
For any sequence $\left(W_n\right)$ a sequence of probability graphons and $W$ a probability graphon, the following are equivalent:
\begin{enumerate}

\item  The sequence $\left(W_n\right)$ is convergent to $W$ in the cut distance $\delta_{\square,\mathcal{LP}};$
\item  The sequence $\left(W_n\right)$ is convergent to $W$ in the cut distance $\delta_{\square,\cF}$ for any choice of convergence determining sequence $\cF;$
\item The overlay functional values $\mathcal{C}\left(W_n, U\right)$ converge to $\mathcal{C}\left(W, U\right)$ for every (Bochner measurable) $\CbFunct-$valued kernel $U$;
\item The overlay functional values $\mathcal{C}\left(W_n, H^{\beta}\right)$ converge to $\mathcal{C}\left(W, H^{\beta}\right)$ for every decorated $\CbFunct-$graph $H^{\beta}$;
\item The overlay functional values $\mathcal{C}\left(W_n, U\right)$ converge to $\mathcal{C}\left(W, U\right)$ for every $\cF-$valued kernel $U$;
\item The overlay functional values $\mathcal{C}\left(W_n, H^{\beta}\right)$ converge to $\mathcal{C}\left(W, H^{\beta}\right)$ for every decorated $\cF-$graph $H^{\beta}$;
\item The quotient sets $\mathcal{Q}_k\left(W_n\right)$ converge $\mathcal{Q}_k\left(W\right)$ to  the $d_{\square}^{\text {Haus }}$ Hausdorff metric for every $k \geq 1$;
\item  The sequence $\left(W_n\right)$ is tight and the $\cF-$overlay functional values $\mathcal{C_{\cF}}\left(W_n, Q\right)=\mathcal{C_{\cF}}\left(Q,W_n\right)$ are convergent to $\mathcal{C_{\cF}}\left(W, Q\right)=\mathcal{C_{\cF}}\left(Q,W\right)$ for every probability graphon $Q.$
\end{enumerate}
\end{theorem}
\proof 
(1) if and only if (2) follows from Corollary 5.6 in \cite{abraham2023probabilitygraphons} (see Theorem \ref{cor:equiv-topo}). \newline
(2) $\Rightarrow$ (3) by Corollary \ref{CorSequencesEqOv}. \newline
(3) $\Rightarrow$ (4) is trivial by identity \eqref{eq:IdentOverlayGraphGraphon}.\newline
(3) $\Rightarrow$ (5) is trivial.\newline
(4) $\Rightarrow$ (6) is trivial.\newline
(5) $\Rightarrow$ (6) is by identity \eqref{eq:IdentOverlayGraphGraphon}.\newline

(6) $\Rightarrow$ (1) 
We prove the contrapositive. Thus we assume that the sequence of probability graphons $W_n$ is not convergent to the probability graphon $W$ in $\delta_{\square,\mathcal{LP}}.$ This can happen in the following distinct cases:
\begin{itemize}
    \item \emph{Case 1:} The sequence of probability graphons $\left(W_n\right)$ converges to a probability graphon $V$ that is not weakly isomorphic to $W.$
    \item \emph{Case 2:} The sequence of probability graphons $\left(W_n\right)$ is not tight.
    \item \emph{Case 3:} The sequence of probability graphons $\left(W_n\right)$ is tight but not convergent.
\end{itemize}
We will show that the real-valued sequence $\cC(W_n,U)$ does not converge to $\cC(W,U)$ in each case. \newline

\emph{Case 1:} 
Let $\left(W_n\right)$ be a sequence of probability graphons which converges to a probability graphon $V$ that is not weakly isomorphic to $W.$ By Lemma \ref{LemmmetricConvImpliesOverlay}, we have that $$
\cC(W_n,H^{\beta})\rightarrow \cC(V,H^{\beta})
$$
for every $\cF-$graph $H^{\beta}.$ However, by Lemma \ref{EquivWeakIsoOverl} there exists a $\cF-$graph $H^{\beta}$ such that 
$$\cC(V,H^{\beta})\neq \cC(W,H^{\beta}).$$ This concludes the proof in this case. \newline

\emph{Case 2:}
Let $\left(W_n\right)$ be a sequence of probability graphons which is not tight. Therefore, by Theorem \ref{ThmRelCompact}, there does not exist a probability measure $\mu$ such that $\mu_n=M_{W_n}=M_{W^{\varphi_n}_n}$ weakly converges to $\mu.$ 

This means that for every probability measure $\mu$ there exists an $f\in \cF$ such that the real-valued sequence $$
\int_{[0,1]^2}W_n[f]=M_{W_n}[f]=\mu_n[f]=\int_{\Space}f\mathrm{d}\mu_n
$$
does not converge to $\mu[f].$ 

For $f\in \cF$ we denote by $U_f$ the $\cF-$valued kernel $U_f=f\mathbbm{1}_{[0,1]^2}.$

We obtain that there does not exist a probability graphon $W$ such that 
$$
\cC(W_n,U_f)=\sup_{\varphi}\int_{[0,1]^2}W_n(U_f^{\varphi})=\int_{[0,1]^2}W_n[f]=M_{W_n}[f]\rightarrow M_W[f]=\cC(W,U_f)$$
 for every $f\in \cF$ and this concludes the proof in this case. \newline

\emph{Case 3:}
The proof is the same as the proof of implication (6) $\Rightarrow$ (1) of  Theorem \ref{ThmEquivalenceConvQuotientOverlay}.\newline

This already proves that (1),(2),(3),(4),(5),(6).\newline

(1) $\Rightarrow$ (7) by Lemma \ref{LemmMetricconvImplQuotient}. \newline
(7) $\Rightarrow$ (6) The proof is the same as the implication (7) $\Rightarrow$ (6) in Theorem \ref{ThmEquivalenceConvQuotientOverlay} where we substitute $W_m$ with $W.$\newline

Therefore, we obtain that (1), (2), (3), (4), (5), (6), (7) are equivalent.\newline
(3) $\Rightarrow$ (8) By Lemma \ref{Eq:FOverlayCOverlayInfty}. \newline
(8) $\Rightarrow$ (1) 
We prove the contrapositive. Thus we assume that the sequence of probability graphons $W_n$ is not convergent to the probability graphon $W$ in $\delta_{\square,\mathcal{LP}}.$ This can happen in the following distinct cases:
\begin{itemize}
    \item \emph{Case 1:} The sequence of probability graphons $\left(W_n\right)$ converges to a probability graphon $V$ that is not weakly isomorphic to $W.$
    \item \emph{Case 2:} The sequence of probability graphons $\left(W_n\right)$ is not tight.
    \item \emph{Case 3:} The sequence of probability graphons $\left(W_n\right)$ is tight but not convergent.
\end{itemize}
We will show that the real-valued sequence $\cC_{\cF}(W_n,U)$ does not converge to $\cC(W,U)$ in each case. \newline

\emph{Case 1:} 
Let $\left(W_n\right)$ be a sequence of probability graphons which converges to a probability graphon $V$ that is not weakly isomorphic to $W.$ By Lemma \ref{LemmmetricConvImpliesOverlay}, we have that $$
\cC_{\cF}(W_n,Q)\rightarrow \cC_{\cF}(V,Q)
$$
for every probability graphon $Q.$ However, by Lemma \ref{EquivWeakIsoOverl} there exists a probability graphon $Q$ such that 
$$\cC_{\cF}(V,Q)\neq \cC_{\cF}(W,Q).$$ This concludes the proof in this case. \newline

\emph{Case 2:}
Let $\left(W_n\right)$ be a sequence of probability graphons which is not tight. Therefore, by Theorem \ref{ThmRelCompact}, there does not exist a probability measure $\mu$ such that $\mu_n=M_{W_n}=M_{W^{\varphi_n}_n}$ weakly converges to $\mu.$ 

This means that for every probability measure $\mu$ there exists an $f\in \cF=(f)_{k\in \N}$ such that the real-valued sequence $$
\int_{[0,1]^2}W_n[f]=M_{W_n}[f]=\mu_n[f]=\int_{\Space}f\mathrm{d}\mu_n
$$
does not converge to $\mu[f].$ 

For a probability measure $\nu \in \Proba,$ we denote by $U_{\nu}$ the probability graphon  $U_{\nu}=\nu\mathbbm{1}_{[0,1]^2}.$

We obtain that 
$$
\cC_{\cF}(W_n,U_{\nu})=\sup_{\varphi}\sum^{\infty}_{k=0}2^{-k}\int_{[0,1]^2}W_n[f_k]U_{\nu}^{\varphi}[f_k]=\sum^{\infty}_{k=0}2^{-k}\nu[f_k]\int_{[0,1]^2}W_n[f_k]=\sum^{\infty}_{k=0}2^{-k}\nu[f_k]M_{W_n}[f_k]$$
 for every $\nu\in \Proba.$  

Let's now suppose that there exists a probability graphon $W$ such that 

$$
\cC_{\cF}(W_n,U)-\cC_{\cF}(W,U)\rightarrow 0.
$$
Therefore, by the monotone convergence theorem, we should have 
$$
\begin{aligned}
&    \left(\cC_{\cF}(W_n,U_{\nu})-\cC_{\cF}(W,U_{\nu})\right)=\sum^{\infty}_{k=0}2^{-k}\nu[f_k]\left(M_{W_n}[f_k]-M_{W}[f_k]\right)=\\
&\sum^{\infty}_{k=0}2^{-k}\nu[f_k]\left(M_{W_n}-M_{W}\right)[f_k]=\int_{\Space}\sum^{\infty}_{k=0}\space2^{-k}\left(M_{W_n}-M_{W}\right)[f_k]f_k(z) \nu(dz)
\end{aligned}
$$
Therefore, choosing $\nu=\delta_z$ for every $z\in \Space,$ we obtain the sequence of continuous and bounded functions 
$$
s_n(z)=\sum^{\infty}_{k=0}\space2^{-k}\left(M_{W_n}-M_{W}\right)[f_k]f_k(z) \rightarrow0.
$$

However, applying $\left(M_{W_n}-M_{W}\right)$ to $s_n$ and using the dominated convergence theorem, we obtain  
$$
\left(M_{W_n}-M_{W}\right)[s_n]=\sum^{\infty}_{k=0}\space2^{-k}\left(\left(M_{W_n}-M_{W}\right)[f_k]\right)^2\rightarrow 0.
$$
This implies $M_{W_n}[f]\rightarrow M_{W}[f] $ for every $f\in \cF=(f_k)_{k\in \N}.$ But this is in contradiction with the initial assumption that for every measure $\mu$ there exists an $f\in \cF$ such that $M_{W_n}[f]$ does not converge to $\mu[f].$  \newline
\emph{Case 3:}
The proof is the same as the proof of implication (8) $\Rightarrow$ (1) of  Theorem \ref{ThmEquivalenceConvQuotientOverlay}.\newline

\endproof

\begin{remark}
Similar considerations as Remark \ref{RemarkMainThmMorGen} apply to Theorem \ref{ThmEquivalenceConvQuotientOverlayVersion2}.
\end{remark}

\textbf{Acknowledgements:}  The author would like to thank Julien Weibel for clarifying results about probability graphons.
\medskip
\section*{References}

\bibliographystyle{plain} 
\bibliography{biblio}

\begin{thebibliography}{10}

\bibitem{abraham2023probabilitygraphons}
R.~Abraham, J-F. Delmas, and J.~Weibel.
\newblock Probability-graphons: Limits of large dense weighted graphs, 2023. arXiv:2312.15935v1 [cs.DM].

\bibitem{aldous1981representations}
D.J. Aldous.
\newblock Representations for partially exchangeable arrays of random variables.
\newblock {\em Journal of Multivariate Analysis}, 11(4):581--598, 1981.

\bibitem{aldous2010exchangeability}
D.J. Aldous.
\newblock Exchangeability and continuum limits of discrete random structures.
\newblock In {\em Proceedings of the International Congress of Mathematicians 2010 (ICM 2010) (In 4 Volumes) Vol. I: Plenary Lectures and Ceremonies Vols. II--IV: Invited Lectures}, pages 141--153. World Scientific, 2010.

\bibitem{ayi2023graphlimitinteractingparticle}
N.~Ayi and N.~Pouradier Duteil.
\newblock Graph limit for interacting particle systems on weighted random graphs, 2023.

\bibitem{ayi2024meanfieldlimitnonexchangeablemultiagent}
N.~Ayi, N.~Pouradier Duteil, and D.~Poyato.
\newblock Mean-field limit of non-exchangeable multi-agent systems over hypergraphs with unbounded rank, 2024.

\bibitem{backhausz2018action}
Á. Backhausz and B.~Szegedy.
\newblock Action convergence of operators and graphs.
\newblock {\em Canadian Journal of Mathematics}, 74(1):72–121, 2022.

\bibitem{barthelemy_2016_structure}
M.~Barthelemy.
\newblock {\em The structure and dynamics of cities: {U}rban data analysis and theoretical modeling}.
\newblock Cambridge University Press, Cambridge, UK, 2016.

\bibitem{HypergraphsNetworks}
F.~Battiston, G.~Cencetti, I.~Iacopini, V.~Latora, M.~Lucas, A.~Patania, J.~Young, and G.~Petri.
\newblock Networks beyond pairwise interactions: Structure and dynamics.
\newblock {\em Physics Reports}, 874:1--92, 2020.

\bibitem{HigerOrdIntBook}
F.~Battiston and G.~Petri.
\newblock {\em Higher-Order Systems}.
\newblock Understanding Complex Systems. Springer Cham, 2022.

\bibitem{battiston_2012_debtrank}
S.~Battiston, M.~Puliga, R.~Kaushik, P.~Tasca, and G.~Caldarelli.
\newblock {DebtRank}: {T}oo central to fail? {F}inancial networks, the {FED} and systemic risk.
\newblock {\em Scientific Reports}, 2(1):541--541, 2012.

\bibitem{BenjaminiLimit}
I.~Benjamini and O.~Schramm.
\newblock Recurrence of distributional limits of finite planar graphs.
\newblock {\em Electronic Journal of Probability}, 6:1 -- 13, 2001.

\bibitem{bick2024dynamicalGraphLimi}
C.~Bick and D.~Sclosa.
\newblock Dynamical systems on graph limits and their symmetries.
\newblock {\em Journal of Dynamics and Differential Equations}, pages 1--36, 2024.

\bibitem{billingsley1968convergence}
P.~Billingsley.
\newblock {\em Convergence of Probability Measures}.
\newblock Wiley Series in Probability and Mathematical Statistics. Wiley, 1968.

\bibitem{BogachevMT2}
V.I. Bogachev.
\newblock {\em Measure Theory}.
\newblock Number v. 2 in Measure Theory. Springer, 2006.

\bibitem{BogachevMT1}
V.I. Bogachev.
\newblock {\em Measure Theory}.
\newblock Number v. 1 in Measure Theory. Springer Berlin Heidelberg, 2007.

\bibitem{Bogachev}
V.I. Bogachev.
\newblock {\em Weak Convergence of Measures}.
\newblock Mathematical Surveys and Monographs. American Mathematical Society, 2018.

\bibitem{local-global1}
B.~Bollob\'{a}s and O.~Riordan.
\newblock Sparse graphs: Metrics and random models.
\newblock {\em Random Struct. Algorithms}, 39(1):1–38, 2011.

\bibitem{BorgsNonparaStat}
C.~Borgs and J.~Chayes.
\newblock Graphons: A nonparametric method to model, estimate, and design algorithms for massive networks.
\newblock In {\em Proceedings of the 2017 ACM Conference on Economics and Computation}, EC '17, page 665–672, New York, NY, USA, 2017. Association for Computing Machinery.

\bibitem{borgs2011convergentAnnals}
C.~Borgs, J.~Chayes, L.~Lovasz, V.T. Sos, and K.~Vesztergombi.
\newblock Convergent sequences of dense graphs {II}: Multiway cuts and statistical physics.
\newblock {\em Annals of Mathematics}, December 2011.

\bibitem{BORGS20081801}
C.~Borgs, J.~Chayes, L.~Lovász, V.T. Sós, and K.~Vesztergombi.
\newblock Convergent sequences of dense graphs {I: Subgraph} frequencies, metric properties and testing.
\newblock {\em Advances in Mathematics}, 219(6):1801--1851, 2008.

\bibitem{bramburger2023pattern}
J.~Bramburger and M.~Holzer.
\newblock Pattern formation in random networks using graphons.
\newblock {\em SIAM Journal on Mathematical Analysis}, 55(3):2150--2185, 2023.

\bibitem{SubModLovaszquotientconvergence}
K.~Bérczi, M.~Borbényi, L.~Lovász, and L.M. Tóth.
\newblock Quotient-convergence of submodular setfunctions, 2024.

\bibitem{Bohle_2021}
T.~Böhle, C.~Kuehn, R.~Mulas, and J.~Jost.
\newblock Coupled hypergraph maps and chaotic cluster synchronization.
\newblock {\em Europhysics Letters}, 136(4):40005, mar 2022.

\bibitem{carletti2020dynamical}
T.~Carletti, D.~Fanelli, and S.~Nicoletti.
\newblock Dynamical systems on hypergraphs.
\newblock {\em Journal of Physics: Complexity}, 1(3):035006, 2020.

\bibitem{MeasurableSetsMeasures}
L.E. Dubins and D.~A. Freedman.
\newblock Measurable sets of measures.
\newblock {\em Pacific Journal of Mathematics}, 14:1211--1222, 1964.

\bibitem{StrongMeasEquivMeasSeparNBanach}
L.E. Dubins and D.A. Freedman.
\newblock Probability theory in banach spaces: An introductory survey.
\newblock In A.T. Bharucha-Reid, editor, {\em Random Integral Equations}, volume~96 of {\em Mathematics in Science and Engineering}, pages 7--63. Elsevier, 1972.

\bibitem{HypergraphsSzegedy2}
G.~Elek and B.~Szegedy.
\newblock Limits of hypergraphs, removal and regularity lemmas. a non-standard approach, 2007. arXiv:0705.2179 [math.CO].

\bibitem{hypergrELEK20121731}
G.~Elek and B.~Szegedy.
\newblock A measure-theoretic approach to the theory of dense hypergraphs.
\newblock {\em Advances in Mathematics}, 231(3):1731--1772, 2012.

\bibitem{Ethier}
S.~N. Ethier and T.~G. Kurtz.
\newblock {\em Markov processes -- characterization and convergence}.
\newblock Wiley Series in Probability and Mathematical Statistics: Probability and Mathematical Statistics. John Wiley \& Sons Inc., New York, 1986.

\bibitem{falgasravry2016multicolour}
V.~Falgas-Ravry, K.~O'Connell, J.~Strömberg, and A.~Uzzell.
\newblock Multicolour containers and the entropy of decorated graph limits, 2016.

\bibitem{fornito_2016_fundamentals}
A.~Fornito, A.~Zalesky, and E.~T. Bullmore.
\newblock {\em Fundamentals of brain network analysis}.
\newblock Academic Press, London, UK, first edition, 2016.

\bibitem{frenkel2018convergence}
P.~Frenkel.
\newblock Convergence of graphs with intermediate density.
\newblock {\em Transactions of the American Mathematical Society}, 370(5):3363--3404, 2018.

\bibitem{ComparingProbMeas}
A.~L. Gibbs and F.~E. Su.
\newblock On choosing and bounding probability metrics.
\newblock {\em International Statistical Review}, 70(3):419--435, 2002.

\bibitem{gkogkas2022mean}
M.~A. Gkogkas, C.~Kuehn, and C.~Xu.
\newblock Mean field limits of co-evolutionary heterogeneous networks.
\newblock {\em arXiv preprint arXiv:2202.01742}, 2022.

\bibitem{Hatami2014LimitsOL}
H.~Hatami, L.~Lov{\'a}sz, and B.~Szegedy.
\newblock Limits of locally–globally convergent graph sequences.
\newblock {\em Geometric and Functional Analysis}, 24:269--296, 2014.

\bibitem{hausmann_2013_atlas}
R.~Hausmann, C.~Hidalgo, S.~Bustos, M.~Coscia, S.~Chung, J.~Jimenez, A.~Simoes, and M.~Yildirim.
\newblock {\em The Atlas of Economic Complexity}.
\newblock Puritan Press, 2011.

\bibitem{hoover1979relations}
D.~N. Hoover.
\newblock Relations on probability spaces and arrays of random variables.
\newblock {\em Institute for Advanced Study}, 1979.

\bibitem{ArankaAction2022}
A.~Hrušková.
\newblock Limits of action convergent graph sequences with unbounded $(p,q)$-norms, 2022, arXiv:2210.10720 [math.CO].

\bibitem{jabin2021meanfield}
P.E. Jabin, D.~Poyato, and J.~Soler.
\newblock Mean-field limit of non-exchangeable systems, 2021. arXiv:2112.15406 [math.PR].

\bibitem{jansonGraphonsCutNorm2013}
S.~Janson.
\newblock Graphons, cut norm and distance, couplings and rearrangements.
\newblock {\em New York Journal of Mathematics}, 2013.

\bibitem{JOST2019870}
J.~Jost and R.~Mulas.
\newblock Hypergraph laplace operators for chemical reaction networks.
\newblock {\em Advances in Mathematics}, 351:870--896, 2019.

\bibitem{JostMulasBook}
J.~Jost, R.~Mulas, and D.~Zhang.
\newblock {\em Spectra of Discrete Structures}.
\newblock Under review, 2023.

\bibitem{MedvedevGraphLimitDynamics2}
D.~Kaliuzhnyi-Verbovetskyi and G.~S. Medvedev.
\newblock The mean field equation for the kuramoto model on graph sequences with non-lipschitz limit.
\newblock {\em SIAM Journal on Mathematical Analysis}, 50(3):2441--2465, 2018.

\bibitem{keriven2024functions}
N.~Keriven and S.~Vaiter.
\newblock What functions can graph neural networks compute on random graphs? the role of positional encoding.
\newblock {\em Advances in Neural Information Processing Systems}, 36, 2024.

\bibitem{Kuehn_2020GraphlimitDynam1}
C.~Kuehn.
\newblock Network dynamics on graphops.
\newblock {\em New Journal of Physics}, 22(5):053030, may 2020.

\bibitem{Kuehn_2019GraphlimitDynam2}
C.~Kuehn and S.~Throm.
\newblock Power network dynamics on graphons.
\newblock {\em SIAM Journal on Applied Mathematics}, 79(4):1271--1292, 2019.

\bibitem{kuehn2022vlasov}
C.~Kuehn and C.~Xu.
\newblock Vlasov equations on digraph measures.
\newblock {\em Journal of Differential Equations}, 339:261--349, 2022.

\bibitem{Kunszenti_Kov_cs_2022}
D.~Kunszenti-Kov{\'{a}}cs, L.~Lov{\'{a}}sz, and B.~Szegedy.
\newblock Multigraph limits, unbounded kernels, and banach space decorated graphs.
\newblock {\em Journal of Functional Analysis}, 282(2):109284, jan 2022.

\bibitem{KUNSZENTIKOVACS20191}
D.~Kunszenti-Kovács, L.~Lovász, and B.~Szegedy.
\newblock Measures on the square as sparse graph limits.
\newblock {\em Journal of Combinatorial Theory, Series B}, 138:1--40, 2019.

\bibitem{KUNSZENTIKOVACS2022109284}
D.~Kunszenti-Kovács, L.~Lovász, and B.~Szegedy.
\newblock Multigraph limits, unbounded kernels, and banach space decorated graphs.
\newblock {\em Journal of Functional Analysis}, 282(2):109284, 2022.

\bibitem{MarkovSpaces}
D.~Kunszenti-Kovács, L.~Lovász, and B.~Szegedy.
\newblock Subgraph densities in markov spaces.
\newblock {\em Advances in Mathematics}, 437:109414, 2024.

\bibitem{NEURIPS2023_8154c89c}
T.~Le and S.~Jegelka.
\newblock Limits, approximation and size transferability for gnns on sparse graphs via graphops.
\newblock In A.~Oh, T.~Naumann, A.~Globerson, K.~Saenko, M.~Hardt, and S.~Levine, editors, {\em Advances in Neural Information Processing Systems}, volume~36, pages 41305--41342. Curran Associates, Inc., 2023.

\bibitem{le2023poincar}
T.~Le, L.~Ruiz, and S.~Jegelka.
\newblock A poincar$\backslash$'e inequality and consistency results for signal sampling on large graphs.
\newblock {\em arXiv preprint arXiv:2311.10610}, 2023.

\bibitem{levie2024graphon}
R.~Levie.
\newblock A graphon-signal analysis of graph neural networks.
\newblock {\em Advances in Neural Information Processing Systems}, 36, 2024.

\bibitem{Lovsz2007SzemerdisLF}
L.~Lov{\'a}sz and B.~Szegedy.
\newblock Szemer{\'e}di’s lemma for the analyst.
\newblock {\em GAFA Geometric And Functional Analysis}, 17:252--270, 2007.

\bibitem{LovaszGraphLimits}
L.~Lovász.
\newblock {\em Large Networks and Graph Limits.}, volume~60 of {\em Colloquium Publications}.
\newblock American Mathematical Society, 2012.

\bibitem{LOVASZ2006933}
L.~Lovász and B.~Szegedy.
\newblock Limits of dense graph sequences.
\newblock {\em Journal of Combinatorial Theory, Series B}, 96(6):933--957, 2006.

\bibitem{lovász2010limits}
L.~Lovász and B.~Szegedy.
\newblock Limits of compact decorated graphs, 2010. arXiv:1010.5155 [math.CO].

\bibitem{majhi2022dynamics}
S.~Majhi, M.~Perc, and D.~Ghosh.
\newblock Dynamics on higher-order networks: A review.
\newblock {\em Journal of the Royal Society Interface}, 19(188):20220043, 2022.

\bibitem{RandomMatrixGraphonmale2014}
C.~Male and S.~P{\'e}ch{\'e}.
\newblock Uniform regular weighted graphs with large degree: Wigner's law, asymptotic freeness and graphons limit.
\newblock {\em arXiv preprint arXiv:1410.8126}, 2014.

\bibitem{maskey2023transferability}
S.~Maskey, R.~Levie, and G.~Kutyniok.
\newblock Transferability of graph neural networks: an extended graphon approach.
\newblock {\em Applied and Computational Harmonic Analysis}, 63:48--83, 2023.

\bibitem{maskey2022generalization}
S.~Maskey, R.~Levie, Y.~Lee, and G.~Kutyniok.
\newblock Generalization analysis of message passing neural networks on large random graphs.
\newblock {\em Advances in neural information processing systems}, 35:4805--4817, 2022.

\bibitem{MedvedevGraphLimitDynamics1}
G.~S. Medvedev.
\newblock The nonlinear heat equation on dense graphs and graph limits.
\newblock {\em SIAM Journal on Mathematical Analysis}, 46(4):2743--2766, 2014.

\bibitem{MHJ}
R.~Mulas, D.~Horak, and J.~Jost.
\newblock {Graphs, simplicial complexes and hypergraphs: Spectral theory and topology}.
\newblock In F.~Battiston and G.~Petri, editors, {\em Higher order systems}. Springer, 2022.

\bibitem{HypergraphsDynamics}
R.~Mulas, C.~Kuehn, and J.~Jost.
\newblock Coupled dynamics on hypergraphs: Master stability of steady states and synchronization.
\newblock {\em Phys. Rev. E}, 101:062313, 2020.

\bibitem{MeasTheorActionZucal}
R.~Mulas and G.~Zucal.
\newblock A measure-theoretic representation of graphs.
\newblock {\em Periodica Mathematica Hungarica}, 88:8--24, 2024.

\bibitem{pagani_2013_power}
G.~A. Pagani and M.~Aiello.
\newblock The power grid as a complex network: {A} survey.
\newblock {\em Physica A: Statistical Mechanics and its Applications}, 392(11):2688--2700, 2013.

\bibitem{pastor-satorras_2015_epidemic}
R.~Pastor-Satorras, C.~Castellano, P.~Van~Mieghem, and A.~Vespignani.
\newblock Epidemic processes in complex networks.
\newblock {\em Rev. Mod. Phys.}, 87:925--979, 2015.

\bibitem{Varadarajan}
V.S. Varadarajan.
\newblock Weak convergence of measures on separable metric spaces.
\newblock {\em Sankhyā: The Indian Journal of Statistics (1933-1960)}, 19(1/2):15--22, 1958.

\bibitem{wolfe2013nonparametricgraphonestimation}
P.~J. Wolfe and S.~C. Olhede.
\newblock Nonparametric graphon estimation, 2013.

\bibitem{HypergraphonsZhao}
Y.~Zhao.
\newblock Hypergraph limits: A regularity approach.
\newblock {\em Random Structures and Algorithms}, 47, 03 2014.

\bibitem{RandomMatricesGraphonsZhu}
Y.~Zhu.
\newblock A graphon approach to limiting spectral distributions of wigner-type matrices.
\newblock {\em Random Structures \& Algorithms}, 56(1):251--279, 2020.

\bibitem{zucal2023action}
G.~Zucal.
\newblock Action convergence of general hypergraphs and tensors, 2023. arXiv:2308.00226 [math.CO].

\bibitem{NewPreprintGiulio}
G.~Zucal.
\newblock Graph convergence and graph limits using probability spaces, In preparation.

\end{thebibliography}
\end{document}